\documentclass[12pt,a4paper]{article}
\usepackage{latexsym,epsfig,fullpage,amsfonts,amssymb,amsmath,amscd,epic,graphics,
rotating,theorem,yhmath,epstopdf}
\usepackage{latexsym,fullpage,amsfonts,amssymb,amsmath,amscd,graphics,epic,theorem}
\usepackage[all]{xy}
\usepackage[usenames]{color}
\baselineskip 18pt \textwidth 15cm  \theoremstyle{plain}
\newtheorem{lemma}{Lemma}%[section]
\newtheorem{proposition}[lemma]{Proposition}\newtheorem{conjecture}{Conjecture}
\newtheorem{remark}[lemma]{Remark}
\newtheorem{example}[lemma]{Example}
\newtheorem{theorem}{Theorem}%[section]

{\theorembodyfont{\rmfamily}
\begin{document}
\newcommand{\Int}{{\operatorname{Int}}}\newcommand{\MCTI}{{\operatorname{MCTI}}}
\newcommand{\AG}{{\mathrm A}\Gamma}\newcommand{\MCTE}{{\operatorname{MCTE}}}
\newcommand{\Ima}{{\operatorname{Im}}}
\newcommand{\Rea}{{\operatorname{Re}}}
\newcommand{\BH}{{\operatorname{bh}}}
\newcommand{\GW}{{\operatorname{GW}}}
\newcommand{\ini}{{\operatorname{ini}}}
\newcommand{\fm}{{\mathfrak{m}}}
\newcommand{\Tor}{{\operatorname{Tor}}}
\newcommand{\conv}{{\operatorname{Conv}}}
\newcommand{\Conj}{{\operatorname{conj}}}
\newcommand{\Ev}{{\operatorname{Ev}}}
\newcommand{\pr}{{\operatorname{pr}}}
\newcommand{\val}{{\operatorname{val}}}
\newcommand{\eps}{{\varepsilon}}
\newcommand{\DD}{\boldsymbol{D}}
\newcommand{\idim}{{\operatorname{idim}}}
\newcommand{\DP}{{\operatorname{DP}}}
\newcommand{\Arc}{{\operatorname{Arc}}}
\newcommand{\Graph}{{\operatorname{Graph}}}
\newcommand{\sm}{{\mathrm{sm}}}
\newcommand{\Z}{{\mathbb Z}}
\newcommand{\bw}{{\boldsymbol{w}}}
\newcommand{\R}{{\mathbb{R}}}\newcommand{\D}{{\mathbb{D}}}
\newcommand{\bz}{{\boldsymbol{z}}}
\newcommand{\bp}{{\boldsymbol{p}}}
\newcommand{\Hom}{{\operatorname{Hom}}}
\newcommand{\const}{{\operatorname{const}}}
\newcommand{\Tors}{{\operatorname{Tors}}}
\newcommand{\Spec}{{\operatorname{Spec}}}
\newcommand{\Aut}{{\operatorname{Aut}}}
\newcommand{\C}{{\mathbb{C}}}\newcommand{\Q}{{\mathbb{Q}}}
\newcommand{\bA}{{\boldsymbol{A}}}\newcommand{\bF}{{\boldsymbol{F}}}
\newcommand{\bB}{{\boldsymbol{B}}}
\newcommand{\bn}{{\boldsymbol{n}}}
\newcommand{\bs}{{\boldsymbol{s}}}
\newcommand{\PP}{{\mathbb{P}}}
\newcommand{\Id}{{\operatorname{Id}}}
\newcommand{\Sing}{{\operatorname{Sing}}}
\newcommand{\codim}{{\operatorname{codim}}}
\newcommand{\Ann}{{\operatorname{Ann}}}
\newcommand{\ord}{{\operatorname{ord}}}
\newcommand{\mt}{{\operatorname{mt}}}\newcommand{\Br}{{\operatorname{Br}}}
\newcommand{\Prec}{{\operatorname{Prec}}}\newcommand{\ReBr}{{\operatorname{ReBr}}}
\newcommand{\ImBr}{{\operatorname{ImBr}}}
\newcommand{\proofend}{\hfill$\blacksquare$\bigskip}
\newcommand{\pperp}{{\perp{\hskip-0.4cm}\perp}}
\newcommand{\bk}{{\boldsymbol{k}}}
\newcommand{\bl}{{\boldsymbol{l}}}

%\textsf{doublespace}
% Enter full title and short title for running headers
\title{Morsifications of real plane curve singularities}
%\shorttitle{Welschinger-type invariants}

\author{Peter Leviant\footnote{School of Mathematical Sciences, Tel Aviv University, Ramat Aviv, 69978 Tel Aviv, Israel.
E-mail: piterleviant@gmail.com} \and Eugenii Shustin\footnote{
School of Mathematical Sciences, Tel Aviv University, Ramat Aviv, 69978 Tel Aviv, Israel.
E-mail: shustin@post.tau.ac.il}}
\date{}
\maketitle
\medskip

{\it Dedicated to the memory of a great mathematician Egbert Brieskorn}

\begin{abstract}
{\small
A real morsification of a real plane curve singularity is a real deformation given by a family of
real analytic functions having only real Morse critical points with all saddles on the zero level.
We prove the existence of real morsifications for real plane curve singularities having arbitrary real local branches and pairs of complex conjugate branches satisfying some conditions.
This was known before only in the case of all local branches being real (A'Campo, Gusein-Zade).
We also discuss a relation between real morsifications and the topology of singularities,
extending to arbitrary real morsifications the Balke-Kaenders theorem, which states that the A'Campo--Gusein-Zade diagram
associated to a morsification uniquely determines the topological type of a singularity.}
\end{abstract}

%\maketitle

\tableofcontents

\section*{Introduction}

By a {\bf singularity} we always mean a germ $(C,z)\subset\C^2$ of a plane reduced analytic curve at its singular point $z$.
Irreducible components of the germ
$(C,z)$ are called {\bf branches of $(C,z)$}. Let $f(x,y)=0$ be an (analytic) equation of $(C,z)$, where
$f$ is defined in the closed ball $B(z,\eps)\subset\C^2$
of radius $\eps>0$ centered at
$z$. The ball $B(z,\eps)$ is called the {\bf Milnor ball} of $(C,z)$ (and is denoted in the sequel $B_{C,z}$) if
$z$ is the only singular point of $C$ in $B(z,\eps)$, and $\partial B(z,\eta)$ intersects $C$ transversally
for all $0<\eta\le\eps$. A {\bf nodal deformation} of a singularity $(C,z)$ is a
family of analytic curves $C_t=\{f_t(x,y)=0\}$, where $f_t(x,y)$ is analytic in $x,y,t$ for $(x,y)\in B(C,z)$ and
$t$ varying in an open disc $\D_\zeta\subset\C$ of some radius $\zeta>0$ centered at zero, and where $C_0=C$, $C_t$ is smooth along $\partial B_{C,z}$, intersects $\partial
B_{C,z}$ trasversally for all $t\in\D_\zeta$, for any $t\ne0$, the curve $C_t$ has only ordinary nodes in $B_{C,z}$,
and the number of nodes does not depend on $t$.
The maximal number of nodes in a nodal deformation of $(C,z)$ in $B$ equals $\delta(C,z)$, the $\delta$-invariant
(see, for instance, \cite[\S10]{M}).

Let $(C,z)$ be a real singularity, i.e.,
invariant with respect to the complex conjugation, $z\in C$ its real singular point. Denote by
$\ReBr(C,z)$, $\ImBr(C,z)$ the numbers of real branches and the pairs of complex conjugate
branches centered at $z$, respectively. Let $C_t=\{f_t(x,y)=0\}$, $t\in\D_\zeta$, be an equivariant\footnote{Here
and further on, {\it equivariant} means commuting with the complex conjugation.} nodal deformation
of a real singularity $(C,z)$. Its restriction to $t\in[0,\zeta)$ is called
a {\bf real nodal deformation}.
A real nodal deformation is called a {\bf real morsification} of $(C,z)$ if
each function $f_t$, $0<t<\zeta$, has only real critical points in $B(C,z)$, all critical points are Morse, and all the saddle points have the zero critical level. Clearly, then all maxima have positive critical values, and all minima negative ones.

N. A'Campo \cite{AC,AC1,AC2} and S. Gusein-Zade \cite{GZ,GZ1} performed a foundational research on this subject. In particular, they showed
that real morsifications carry a lot of information on singularities and allow one to compute such invariants as the monodromy and intersection form in vanishing homology in a simple and efficient way. However, some questions have remained open, in particular:

\smallskip

{\bf Question:} {\it Does any real plane curve singularity admit a real morsification?}

\smallskip

Our main result is a partial answer to this question.
Before precise formulation, we should mention that an affirmative answer was given before in the case of all branches of $(C,z)$ being real (below
referred to as a {\bf totally real} singularity), see
\cite[Theorem 1]{AC}\footnote{As pointed to us by S. Gusein-Zade, there is a gap in the proof of
\cite[Theorem 1]{AC}: namely, the function in \cite[Formula (1) in page 12]{AC} does not possess
the claimed properties.}
and \cite[Theorem 4]{GZ2} (see also \cite[Section 4.3]{AGV}).
Notice that any topological type of a curve singularity is presented by a totally real singularity, see \cite[Theorem 3]{GZ2}.

Now we give necessary definitions. A singularity is called {\bf Newton non-degenerate}, if in some local coordinates, it is {\bf strictly Newton non-degenerate}, that
is given by an equation $f(x,y)=0$ with a convenient Newton diagram at $z=(0,0)$ and
such that the truncation of $f(x,y)$ to any edge of the
Newton diagram is a quasihomogeneous polynomial without critical points in $(\C^*)^2$
(i.e., it has no multiple binomial factors). We say that a singularity $(C,z)$ is {\bf
admissible
along its tangent line} $L$ if the singularity $(C_L,z)$ formed by the union of the branches of $(C,z)$ tangent to $L$ is as follows:
$(C_L,z)$ is the union of a Newton non-degenerate singularity with a singularity, whose all branches are smooth.

\begin{theorem}\label{t1}
Let $(C,z)$ be a real singularity, ${\mathcal T}(C,z)=\{z_0=z,z_1,...\}$ the vertices of its minimal resolution tree. For any $z_i\in{\mathcal T}(C,z)$ denote by $(C_i,z_i)$ the germ at $z_i$ of the corresponding strict transform of $(C,z)$. If, for any real point $z_i\in{\mathcal T}(C,z)$,
the singularity $(C_i,z_i)$ is admissible
along each of its non-real tangent lines, then
the real singularity $(C,z)$ admits a real morsification.
\end{theorem}

Note that the case of totally real singularities is included, since then the restrictions
asserted in Theorem are empty. We illustrate the range of singularities covered by Theorem \ref{t1} with a few examples.

\begin{example}\label{exmor1}
(1) Any quasihomogeneous (in real coordinates) singularity satisfies the hypotheses of Theorem \ref{t1}, and their morsifications
can be constructed in the same manner as for the totally real singularities even if the singularity contains
complex conjugate branches, see Section \ref{sqh}.

(2) The simplest singularity satisfying the hypotheses of Theorem \ref{t1} and whose morsification is constructed
by a new method suggested in the present paper is a pair of transversal ordinary cuspidal branches, given, for instance, by
an equation $(x^2+y^2)^2+x^5=0$. The real part of its morsification looks as shown in Figure \ref{fmor1}. One can
show that all possible morsifications are isotopic to this one.

\begin{figure}
\setlength{\unitlength}{0.6cm}
\begin{picture}(7,5)(-5,0)
\epsfxsize 60mm \epsfbox{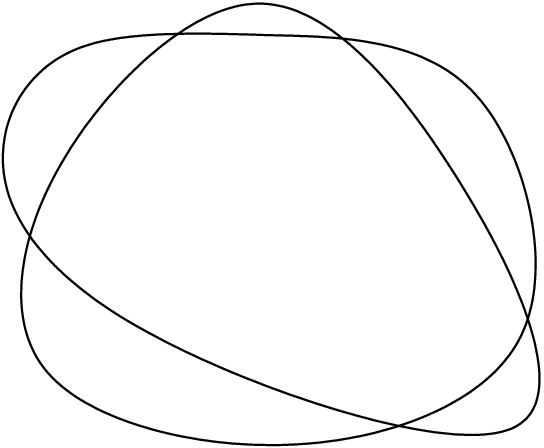}
\end{picture}
\caption{Morsification of a pair of complex conjugate cuspidal branches} \label{fmor1}
\end{figure}

(3) The simplest singularity beyond the range of Theorem \ref{t1} is a pair of two transversal complex conjugate branches
of order $4$ with two Puiseux pairs $(2,3)$ and $(2,7)$ (equivalently, with the Puiseux characteristic
exponents $(4,6,7)$), given, for instance, by an equation
$$((w_+^2-x^3)^2-x^5w_+)((w_-^2-x^3)^2-x^5w_-)=0,\quad w_{\pm}=y\pm x\sqrt{-1}\ .$$
On the other hand, a singularity consisting of a pair of complex conjugate branches with the same Puiseux pairs
$(2,3)$, $(2,7)$ as above, but having a common real tangent does satisfy the hypotheses of Theorem \ref{t1}, since after one blow up
it turns into a singularity with two complex conjugate branches having only one Puiseux pair.
\end{example}

We
believe that the following holds:

\begin{conjecture}\label{con-new}
Any real plane curve singularity possesses a real morsification.
\end{conjecture}

In the proof of Theorem \ref{t1} presented in Section
\ref{sec-exi}, we combine a relatively elementary inductive blow-up
construction in the spirit of \cite{AC} with the patchworking construction as appears in
\cite{Sh1,ST} and some explicit formulas for real morsifications of pairs of complex conjugate
smooth branches and pairs of branches of topological type
$x^p+y^q=0$, $(p,q)=1$. We expect that suitable formulas for real morsifications
of pairs of complex conjugate branches with several Puiseux pairs would lead to a complete solution of the existence problem of real morsifications.

A real morsification
of a totally real singularity yields
a so-called A'Campo-Gusein-Zade diagram, which uniquely determines the
topological type of the singular point, as shown by L. Balke and R. Kaenders
\cite[Theorem 2.5 and Corollary 2.6]{BK}. In Section \ref{sec3}, we extend this result to morsifications of
arbitrary real singularities.

\section{Elementary geometry of real morsifications}\label{sec1}

For the reader's convenience, we present here few simple and in fact known claims on morsifications.
In what follows we consider only real singularities.

Recall that a real node of a real curve can be either hyperbolic or elliptic, that is,
analytically equivalent over $\R$ either to $x^2-y^2=0$, or $x^2+y^2=0$, respectively. For
a real nodal deformation $C_t=\{f_t(x,y)=0\}$, $0\le t<\zeta$, the saddle critical points of
$f_t$ on the zero level correspond to real hyperbolic nodes of $C_t$ and vice versa.

\begin{lemma}\label{l1}
The number of hyperbolic nodes in any real
nodal deformation $C_t$, $0\le t<\zeta$,
of $(C,z)$ does not exceed $\delta(C,z)-\ImBr(C,z)$.
\end{lemma}

{\bf Proof.}
As we noticed in Introduction, the maximal number of nodes in a nodal deformation of a singularity $(C,z)$ is
the $\delta$-invariant $\delta(C,z)$. In a real nodal deformation, a pair $Q,\overline Q$ of complex conjugate branches either glues
up into one surface immersed into $B(C,z)$ thus reducing the total number of nodes by at least one, or $Q$ and $\overline Q$
do not glue up to each other and to other branches and then their intersection points are either complex conjugate nodes or
real elliptic nodes, and, at last, if $Q$ and $\overline Q$ do not glue up to each other, but glue up to some other branches of
$(C,z)$, we loose at least two nodes. So, the bound follows.
\proofend

The following lemma is a version of \cite[Lemma 4 and Theorem 3]{AC}. Let $C_t$, $0\le t<\zeta$, be a a real morsification of a real singularity $(C,z)$.
The sets $\R C_t$, $0<t<\zeta$, are isotopic in the disc $\R B_{C,z}$. Each of them is called a {\bf divide} of the given morsification (more information on divides see in Section \ref{sec-ag}).
Given a divide $D\subset\R B_{C,z}$ of a real morsification of the real singularity $(C,z)$, the connected components of $\R B_{C,z}\setminus D$ disjoint from $\partial\R B_{C,z}$ are called inner components. Denote by $I(D)$ the union of the closures of the inner components of $\R B_{C,z}\setminus D$
(called {\bf body of the divide} in \cite{AC4}).

\begin{lemma}\label{l3}
Let $D=\R C_t$ be a divide of a real morsification of a real singularity $(C,z)$. Then
\begin{enumerate}
\item[(i)] if $(C,z)$ is not a hyperbolic node then $I(D)$ is non-empty, connected, and simply connected;
\item[(ii)] $D$ has $\delta(C,z)-\ImBr(C,z)$ singularities, which are hyperbolic nodes of $C_t$;
\item[(iii)] each inner component of $\R B_{C,z}\setminus D$
is homeomorphic to an open disc;
\item[(iv)] the number $h(C,z)$ of the inner components of $\R B_{C,z}\setminus D$ does not depend on the morsification and satisfies
the relation
$$h(C,z)+\delta(C,z)-\ImBr(C,z)=\mu(C,z)\ ,$$
$\mu(C,z)$ being the Milnor number.
\end{enumerate}
\end{lemma}

{\bf Proof.} In claim (i) suppose that $I(D)$ is not connected.
Then the associated Coxeter-Dynkin
diagram of the singularity $(C,z)$ constructed in \cite{GZ} (see also
\cite[\S3]{GZ1}) appears to be disconnected contrary to the fact that it is always connected \cite{Gab,GZ2}.
Furthermore, $I(D)$ is simply connected since is has no holes by construction.

Statements (ii)-(iv) follow from claim (i), from
the bound $$\#\Sing(D)\le\delta(C,z)-\ImBr(C,z)$$
of Lemma \ref{l1}, from the Milnor formula \cite[Theorem 10.5]{M}
$$\mu(C,z)=2\delta(C,z)-\ReBr(C,z)-2\ImBr(C,z)+1\ ,$$
from the fact that each inner component of $\R B_{C,z}\setminus D$ contains a critical point of
the function $f_t(x,y)$, and hence
$$h(C,z)+\delta(C,z)-\ImBr(C,z)\le\mu(C,z)\ ,$$
and from the calculation of the Euler characteristic of $I(D)$
$$h(C,z)-\left(2\cdot\#\Sing(D)-\ReBr(C,z)\right)+\#\Sing(D)\ge1\ .$$
\proofend

\begin{remark}
In fact, one could equivalently define real morsifications as real nodal deformations having precisely
$\delta(C,z)-\ImBr(C,z)$ hyperbolic nodes as their only singularities.
\end{remark}

\begin{lemma}\label{l2}
Given a real morsification $C_t$, $0\le t<\zeta$, of a real singularity $(C,z)$,
\begin{itemize}\item any real branch $P$ of $(C,z)$ does not glue up with other branches
and deforms into a family of immersed discs $P_t$, $t>0$, whose real point sets $\R P_t\subset\R B_{C,z}$ are immersed segments with
$\delta(P)$ selfintersetions and endpoints on $\partial\R B_{C,z}$;
\item any pair of complex conjugate branches $Q,\overline Q$ of $(C,z)$ do not glue up to other branches, but glue up
to each other so that they deform into a family of immersed cylinders $Q_t$, $t>0$, with the real point set $\R Q_t\subset
\R B_{C,z}$ being an immersed circle disjoint from $\partial B(C,z)$ and having $\delta(Q\cup\overline Q)-1=2\delta(Q)+(Q\cdot\overline Q)-1$ selfintersections (here $(Q\cdot\overline Q)$ denotes the intersection number);
\item for any two real branches $P',P''$, the intersection $\R P'_t\cap\R P''_t$, $t>0$, consists of $(P'\cdot P'')$ points;
\item for any real branch $P$ and a pair of complex conjugate branches $Q,\overline Q$, the intersection
$\R P_t\cap\R Q_t$, $t>0$, consists of $2(P\cdot Q)$ points;
\item for any two pairs of complex conjugate branches $Q',\overline Q'$ and $Q'',\overline Q''$, the intersection
$\R Q'_t\cap\R Q''_t$, $t>0$, consists of $2(Q'\cdot Q'')+2(Q'\cdot\overline Q'')$ points.
\end{itemize}
\end{lemma}

{\bf Proof.}
Straightforward from Lemmas \ref{l1} and \ref{l3}.
\proofend

\begin{lemma}\label{l4}
Let $(C_1,z)$, $(C_2,z)$ be two real singularities without branches in common. If the real singularity
$(C_1\cup C_2,z)$ possesses a real morsification, then each of the real singularities $(C_1,z)$, $(C_2,z)$ possesses a real morsification too.
\end{lemma}

{\bf Proof.}
Straightforward from Lemma \ref{l2}.
\proofend

Given a divide $D$ of a real morsification of a real singularity $(C,z)$,
it follows from Lemma \ref{l3} that
$I(D)$ possesses a
cellular decomposition into $\Sing(D)$ as vertices, the components of $D
\setminus\Sing(D)$, disjoint from $\partial\R B_{C,z}$, as the $1$-cells, and the inner components
of $\R B_{C,z}\setminus D$ as the $2$-cells. Following \cite[\S1]{AC}, we say that the given real morsification defines a {\bf partition}, if, in the above cellular decomposition of $I(D)$, the intersection of the closures of
any two $2$-cells is either empty, or a vertex, or the closure of a $1$-cell.

This property was assumed in the Balke-Kaenders theorem \cite[Theorem 2.5 and Corollary 2.6]{BK} about the recovery
of the topological type of a singularity out of the A'Campo-Gusein-Zade diagram.
In fact, this assumption is not needed (see Section \ref{sec3}).
Here we just notice the following:

\begin{lemma}\label{ap1}
There are real morsifications that do not define a partition.
\end{lemma}

\begin{figure}
\setlength{\unitlength}{0.8cm}
\begin{picture}(14.5,10)(0,-0.7)
\epsfxsize 135mm \epsfbox{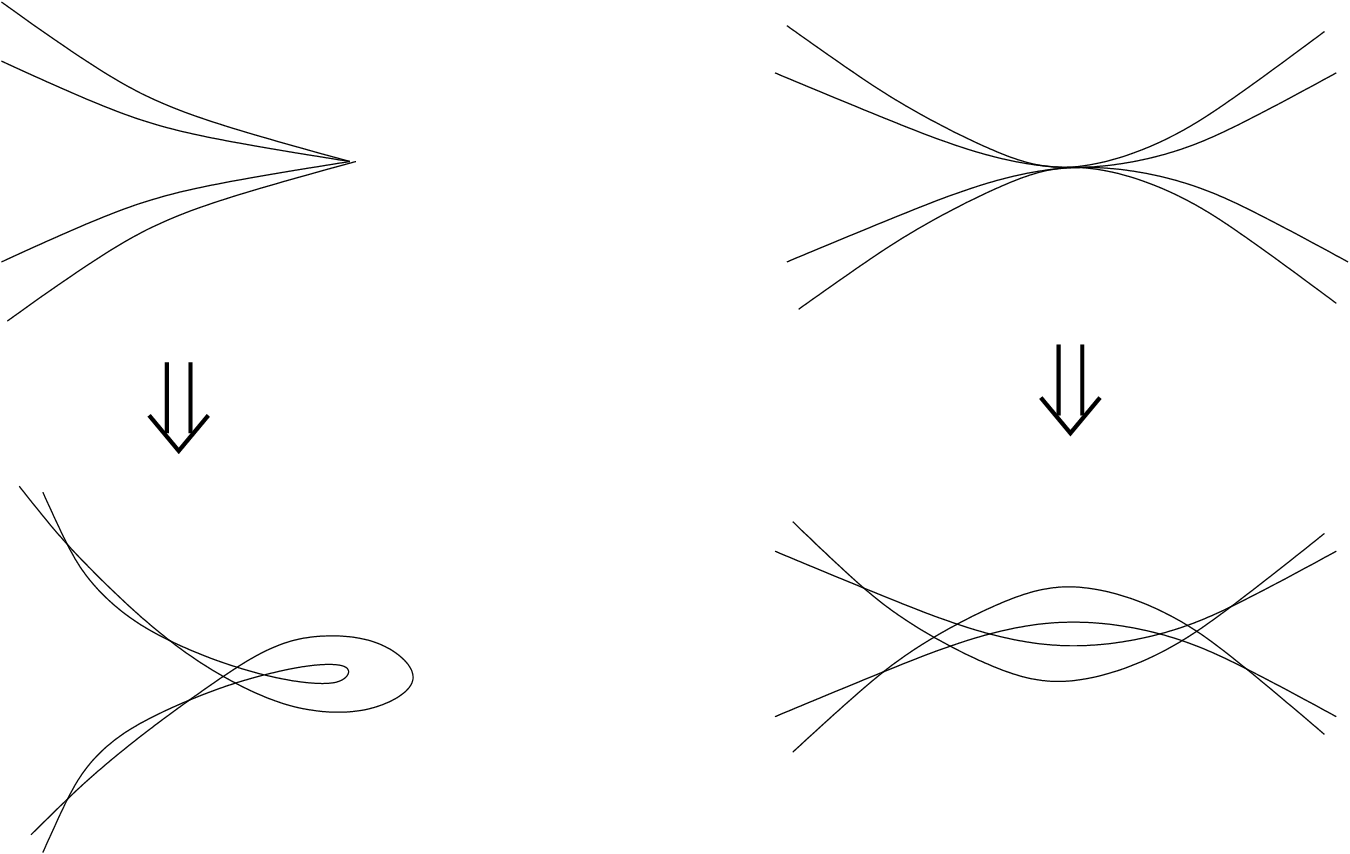}
\put(-14,-0.7){(a)}\put(-4,-0.7){(b)}
\end{picture}
\caption{Non-partitions} \label{fig1}
\end{figure}

{\bf Proof.} For the proof, we present two simple examples: Figure \ref{fig1}(a) shows a real morsification of the
singularity $(y^2+x^3)(y^2+2x^3)=0$ (two cooriented real cuspidal branches with a
common tangent), while Figure \ref{fig1}(b) shows a real morsification of the real
singularity $(y^2-x^4)(y^2-2x^4)=0$ (four real smooth branches quadratically tangent to each other).
A construction is elementary. For example, the morsification shown in Figure \ref{fig1}(a) can be
defined by
$$(y^2+x^2(x-\eps_1(t)))(y^2+2(x-\eps_2(t))^2(x-\eps_3(t)))=0\ ,$$
where $0<\eps_2(t)<\eps_3(t)\ll\eps_1(t)\ll1$.
\proofend

\section{Existence of real morsifications}\label{sec-exi}

\subsection{Blow-up construction}\label{sec2}

Let us recall that the multiplicity of a singularity $(C,z)$, resp. of a branch $P$, is the intersection numbers
$\mt(C,z)=(C\cdot L)_z$, resp. $(P\cdot L)_z$ with a generic line $L$ through $z$. Recall that the proper transform of $(C,z)$ under the blowing up
of $z$ consists of several germs $(C^*_i,z_i)$ with $z_i$ being distinct points on the exceptional divisor $E$
associated with distinct tangents to $(C,z)$. It is know that (see, for instance, \cite[Page 185 and Proposition 3.34]{GLS})
\begin{equation}\delta(C,z)=\sum_i\delta(C^*_i,z_i)+\frac{\mt(C,z)(\mt(C,z)-1)}{2},\quad \mt(C,z)=\sum_i(C^*_i\cdot E)_{z_i}\ .
\label{e2}\end{equation}

\subsubsection{The totally real singularities}\label{total}
The existence of real morsifications for totally real singularities was proved in \cite[Theorem 1]{AC}. We present here a proof (similar to the A'Campo's one)
in order to be self-contained and to use elements of that proof in the general case.

\smallskip

{\bf (1)} Consider, first, the case of a totally real singularity $(C,z)$ whose all branches are smooth.
We proceed by induction on the maximal $\delta$-invariant $\Delta_1(C,z)$ of the union of any subset of branches tangent to each other.

The base of induction, $\Delta_1(C,z)=0$, corresponds to the union of $d\ge2$ smooth branches with distinct tangents.
Here $\delta(C,z)=d(d-1)/2$, and we construct a real morsification by shifting the branches to a general position.

Assuming that $\Delta_1(C,z)>0$ in the induction step, we blow up the point $z$ into an exceptional divisor $E$. The strict transform of $(C,z)$
splits into components $(C^*_i,z_i)$, $z_i\in \R E$, corresponding to different tangents of $(C,z)$.
Notice that $E$ is transversal to all branches of $(C^*_i,z_i)$, and hence
$\Delta_1(C^*_i\cup E,z_i)<\Delta_1(C,z)$ for all $i$ (cf. (\ref{e2})).
Then we construct real morsifications of each real singularity $(C^*_i\cup E,z_i)$ in which the germs
$(E,z_i)$ stay fixed (in view of Lemma \ref{l2} these germs do not glue up with other branches, and hence
can be kept fixed by suitable local equivariant diffeomorphisms). Thus, we get
the union of real curves $(C^*_i)^+$ in neighborhoods of $z_i$, having
$$\sum_i\delta(C^*_i,z_i)=\delta(C,z)-\frac{\mt(C,z)\cdot(\mt(C,z)-1)}{2}$$ real hyperbolic nodes
and $\mt(C,z)$ real intersetion points with $E$. Then we blow down $E$ and obtain a deformation
whose elements have $\delta(C,z)-\frac{\mt(C,z)\cdot(\mt(C,z)-1)}{2}$ real hyperbolic nodes and
a point of transversal intersection of $\mt(C,z)$ smooth branches. Deforming the latter real singularity, we
complete the construction
of a real morsification.

\smallskip

{\bf(2)} Now we prove the existence of real morsifications for arbitrary totally real singularities, using induction
on $\Delta_2(C,z)$, the $\delta$-invariant of the union of all singular branches of $(C,z)$. The preceding
consideration serves as the base of induction. The induction step is very similar: we blow up the point
$z$ and notice that $\sum_i\Delta_2(C^*_i\cup E,z_i)<\Delta_2(C,z)$; then proceed as in the preceding paragraph.

\subsubsection{Semiquasihomogeneous singularities}\label{sqh}
The same blow-up construction of real morsifications works well in
the important particular case of semiquasihomogeneous singularities. Let
$F(x,y)=\sum_{pi+qj=pq}a_{ij}x^iy^j$ be a real square-free quasihomogeneous polynomial, where $1\le p\le q$. Then
$(C,z)=\{F(x,y)+\sum_{pi+qj>pq}a_{ij}x^iy^j=0\}$ is called a real semiquasihomogenelous singularity of type $(p,q)$.
This real singularity has $d=\gcd(p,q)$ branches, among which we allow complex conjugate pairs.

\smallskip

{\bf (1)} A semiquasihomogeneous singularity of type $(p,p)$ is just the union of smooth transversal branches.
If they all are real the existence of a real morsification is proved in Section \ref{total}. Thus, suppose that
$F(x,y)$ splits into the product $F_1(x,y)$ of real linear forms and the product $F_2(x,y)$ of positive definite
quadratic forms $q_i(x,y)$, $1\le i\le k$, $k\ge1$. The forms $q_i$
are not proportional to each other, and there are $b_i>0$, $i=1,...,k$, such that any two quadrics
$q_i-b_i=0$ and $q_j-b_j=0$, $1\le i<j\le k$, intersect in four real points, and all their intersection points are distinct.
So, we obtain a real morsification by deforming $(C,z)$ in the family
$$F(x,y,t)=F_1(x,y)\prod_{i=1}^k(q_i(x,y)-b_it),\quad 0\le t\ll1\ ,$$ and then by shifting each of the lines defined by $F_1=0$ to
a general position.

\smallskip

{\bf (2)} Let $(C,z)$ be a real semiquasihomogeneous singularity of type $(p,q)$, $2\le p<q$.
We simultaneously prove the existence of real morsifications of $(C,z)$ and of the following
additional singularities:
\begin{enumerate}\item[(f1)] $(C\cup L,z)$, where $L$ is a real line intersecting $(C,z)$ at $z$ with multiplicity $p$ (i.e. transversally) or
$q$ (tangent);
\item[(f2)] $(C\cup L_1\cup L_2,z)$, where a real line $L_1$ intersects $(C,z)$ with multiplicity $p$ and a real line
$L_2\ne L_1$ intersects $(C,z)$ at $z$ with multiplicity $p$ or $q$.
\end{enumerate}
We proceed by induction on $\delta(C,z)$. The base of induction, $\delta(C,z)=1$, corresponds to $p=2$, $q=3$, that is, an ordinary cusp.
Here $(C,z)$, $(C\cup L,z)$, and $(C\cup L_1\cup L_2,z)$ are totally real, hence possess a real morsification.
Suppose that $\delta(C,z)>1$, blow up the point $z$, and consider the union of the strict transform of
the studied singularity with the exceptional divisor $E$. Notice that the strict transform of a real semiquasihomogeneous singularity
of type $(p,q)$ is also a real semiquasihomogeneous singularity either of type $(p,q-p)$ if $2p\le q$, or of type
$(q-p,p)$ if $2p>q$, and in both cases it intersects $E$ with multiplicity $p$. It is easy to see that the strict transform
of singularities of the form (f1) and (f2) with added $E$ is again a real singularity of one of these forms
with parameters $(p,q-p)$ or $(q-p,p)$ and, may be, an extra real node. We then complete the proof as in Section \ref{total}.

\subsection{Singularities without real tangents}\label{sec4}
The constructions of morsifications presented in this section is the mein novelty of the present paper. In the case of singularities with only smooth branches,
Lemma \ref{lsmooth} presents a rather simple direct formula for the morsification. In the case of non-smooth branches with one Puiseux pair
(Lemma \ref{lNND} below), we apply an ad hoc deformation argument (a kind of the pathchworking construction). The geometric background
for this argument is as follows. We extend the pair $(\C^2,(C,z))$ to a trivial family $(\C^2,(C,z))\times(\C,0)$, then blow up
the point $z\in\C^2\times\{0\}$. The central fiber of the new family is the union of the blown-up plane
$\C^2_1$ and the exceptional divisor $E\simeq\PP^2$. The germ $(C,z)$ yields in $\PP^2$ a real conic $C_2$ with multiplicity $p\ge2$
that intersects the line $\C^2_1\cap E$ in two imaginary points. Our deformation gives an inscribed
equivariant family of curve germs, whose real part appears to be a deformation of the above $p$-multiple conic $C_2$.

\subsubsection{The case of one pair of complex conjugate tangents}\label{sec-smooth}
Let a real singularity $(C,z)$ have exactly two tangent lines, and they are complex conjugate.
In suitable local equivariant coordinates $x,y$ in $B_{C,z}$, we have $z=(0,0)$, and the tangent lines are
$$L=\{x+(\alpha+\beta\sqrt{-1})y=0\},\
\overline L=\{x+(\alpha-\beta\sqrt{-1})y=0\}\ ,$$
where $\alpha,\beta\in\R$, $\beta\ne0$.

Denote by $(C_i,z)$, $i=1,...,s$, the branches of $(C,z)$ tangent to $L$; respectively $(\overline C_i,z)$, $i=1,...,s$, are the branches of $(C,z)$ tangent to $\overline L$. Introduce the new coordinates $$w=x+(\alpha+\beta\sqrt{-1})y,\quad\widehat w=x+(\alpha-\beta\sqrt{-1})y\ .$$
Notice that $\widehat w=\overline w$ if $x,y\in\R$. We also will use for $\R^2\setminus\{0\}$ the coordinates
$\rho>0$, $\theta\in\R/2\pi\Z$ such that
\begin{equation}x+\alpha y=\rho\cos\theta,\quad \beta y=\rho\sin\theta,\quad
\rho=\sqrt{w\widehat w}\ .\label{e8b}\end{equation}

\begin{lemma}\label{lsmooth}
Let $(C,z)$ have only smooth branches. Then $(C,z)$ possesses a real morsification.
\end{lemma}

{\bf Proof.}
A branch $(C_i,z)$, $1\le i\le s$, has an analytic equation
$$\widehat w=\sum_{n\in I_i}a_{in}w^n,\quad I_i\subset\{n\in\Z\ :\ n>1\},\ a_{in}\in\C^*\ \text{as}\ n\in I_i\ .$$ Correspondingly, $(\overline C_i,z)$ is given by
$w=\sum_{n\in I_i}\overline a_{in}\widehat w^n$. We claim that the equation
\begin{equation}F_t(w,\widehat w):=\prod_{i=1}^s(\Phi_i(w,\widehat w)-t^2)=0,\quad 0\le t<\zeta\ ,
\label{ep1}\end{equation}
defines a real morsification of $(C,z)$, where
$$\Phi_i(w,\widehat w)=\left(\widehat w-\sum_{n\in I_i}a_{in}w^n\right)\left(
w-\sum_{n\in I_i}\overline a_{in}\widehat w^n\right)$$ and $\zeta>0$ is sufficiently small.
First, $F_t(w,\widehat w)$ (the left-hand side of (\ref{ep1}))
is an analytic function in $w,\widehat w$ and $t$.
A separate factor in $F_t(w,\widehat w)$ is
$$\Phi_i(w,\widehat w)-t^2=w\widehat w-t^2+\sum_{n\in I_i}|a_{in}|^2(w\widehat w)^n
-\sum_{n\in I_1}(\overline a_{in}w^{n+1}+a_{in}\widehat w^{n+1})$$
$$+2\sum_{\renewcommand{\arraystretch}{0.6}
\begin{array}{c}
\scriptstyle{n_1<n_2}\\
\scriptstyle{n_1,n_2\in I_i}\end{array}}(w\widehat w)^{n_1}(a_{in_1}\overline a_{in_2}w^{n_2-n_1}+
\overline a_{in_1}a_{in_2}\widehat w^{n_2-n_1})\ .$$
Restricting the equation $\Phi_i(w,\widehat w)$ to $\R B_{C,z}$
(in coordinates $x,y$), passing in $\R^2\setminus\{0\}$ to coordinates $\rho>0$, $\theta\in\R/2\pi\Z$
defined via (\ref{e8b}, and rescaling by substitution of $t\rho$ for $\rho$, we obtain
a family of curves depending on the parameter $0\le t<\zeta$
$$\Psi_{i,t}:=\rho^2-1+\sum_{n\in I_i}t^{2n-2}|a_{in}|^2\rho^{2n}-2\sum_{n\in I_i}t^{n-1}|a_{in}|\rho^{n+1}
\cos((n+1)\theta-\theta_{in})$$
$$+2\sum_{\renewcommand{\arraystretch}{0.6}
\begin{array}{c}
\scriptstyle{n_1<n_2}\\
\scriptstyle{n_1,n_2\in I_i}\end{array}}t^{n_1+n_2-2}|a_{in_1}a_{in_2}|\rho^{n_1+n_2}\cos((n_2-n_1)\theta+
\theta_{in_1}-\theta_{in_2})=0\ ,$$ where $a_{in}=|a_{in}|\exp(\sqrt{-1}\theta_{in})$, $n\in I_i$.
It is easy to see that each of them a circle embedded into
an annulus $\{|\rho-1|<Kt\}\subset\R^2$ with $K>0$ a constant determined by the given singularity $(C,z)$, and, furthermore,
the normal projection of each curve to the circle $\rho=1$ is a diffeomorphism. Let $1\le i<j\le s$. Set $$n_{ij}=\min\{n\in I_i\cup I_j\ :\ a_{in_{ij}}\ne a_{jn_{ij}}\}\ .$$ Note that $n_{ij}=(C_i\cdot C_j)$, the intersection number of branches $C_i,C_j$. On the other hand,
$$\Psi_{i,t}(\rho,\theta)-\Psi_{j,t}(\rho,\theta)=2t^{n_{ij}-1}|a_{in_{ij}}-a_{jn_{ij}}|\rho^{n_{ij}+1}
\cos((n_{ij}+1)\theta-\theta_{ij,n_{ij}})+O(t^{n_{ij}})\ ,$$ where $\theta_{ij,n_{ij}}\in\R/2\pi\Z$, and hence, for
a sufficiently small $t>0$, the curves $\Psi_{i,t}=0$ and $\Psi_{j,t}=0$ intersect transversally in
$2n_{ij}+2$ points. In total, we obtain
$$2\sum_{1\le i< j\le s}(n_{ij}+1)=2\sum_{1\le i<j\le s}(C_i\cdot C_j)+s^2-s
=\delta(C,z)-\ImBr(C,z)$$ hyperbolic nodes as required for a real morsification.
\proofend

\begin{lemma}\label{lNND}
Let the singularity $(C_L,z)$ be formed by a pair of branches of topological type
$x^p+y^q=0$, $2\le p<q$, $(p,q)=1$, that are tangent to $L$ and $\overline L$ respectively.
Then $(C,z)$ possesses a real morsification.
\end{lemma}

{\bf Proof.} {\bf (1)} We start with the very special case of $(C,z)$ given by
\begin{equation}F(w,\widehat w)=w^p\widehat w^p-a\widehat w^{p+q}-\overline a w^{p+q}=0,
\quad a\in\C^*\ .\label{ecrit3}\end{equation}
Denote by $P(\lambda)=\lambda^p+b^{(0)}_{p-2}\lambda^{p-2}+...+b^{(0)}_0
\in\R[\lambda]$ the monic polynomial of degree $p$ having $\left[\frac{p}{2}\right]$ critical
points on the level $-2|a|$ and $\left[\frac{p-1}{2}\right]$ critical points on the level $2|a|$, whose roots sum up to zero (a kind of the $p$-th Chebyshev polynomial). We claim that there exist real functions
$b_0(t),...,b_{p-2}(t)$, analytic in $t^{\frac{1}{p}}$ such that $b_i(0)=b^{(0)}_i$, $0\le i\le p-2$, and the
family
\begin{equation}
F_t(w,\widehat w)=(w\widehat w-t^2)^p+\sum_{i=p-2}^0t^{\frac{(p-i)(p+q)}{p}}b_i(t)(w\widehat w-t^2)^i-a\widehat w^{p+q}-\overline a w^{p+q}=0\ ,
\label{ecrit2}\end{equation} $$0\le t<\zeta\ ,$$ is a real morsification of $(C,z)$. To prove this, we
rescale the latter equation by substituting $(tw,t\widehat w)$ for $(w,\widehat w)$ and restrict our attention
to $\R B_{C,z}$ passing to the coordinates $\rho,\theta$ in (\ref{e8b}):
$$(\rho^2-1)^p+\sum_{i=p-2}^0t^{\frac{(p-i)(q-p)}{p}}b_i(t)(\rho^2-1)^i-2|a|\rho^{p+q}\cos((p+q)\theta-
\theta_a)=0\ ,$$ where $a=|a|\exp(\sqrt{-1}\theta_a)$. Next, we substitute $\rho^2=1+t^{\frac{q-p}{p}}\sigma$ and come to
\begin{equation}
(1+t^{\frac{q-p}{p}}\sigma)^{-(p+q)/2}\left(\sigma^p+\sum_{i=p-2}^0b_i(t)\sigma^i\right)
=2|a|\cos((p+q)\theta-\theta_a)\ .
\label{ecrit1}\end{equation} Finally, we recover the unknown functions $b_{p-2}(t),...,b_0(t)$ from the following conditions.

Let $P(\lambda)>3|a|$ as $|\lambda|>\lambda_0$. Suppose that $|\sigma|\le\lambda_0$ and that $t$ is small so that the function of $\sigma$
$$P_t(\sigma):=(1+t^{\frac{q-p}{p}}\sigma)^{-(p+q)/2}\left(\sigma^p+\sum_{i=p-2}^0b_i(t)\sigma^i\right)$$ has
simple critical points $\mu_1(t),...,\mu_{p-1}(t)$ arranged in the growing order and respectively close to the
critical points $\mu^{(0)}_1,...,\mu^{(0)}_{p-1}$ of $P(\lambda)$. So, we require
\begin{equation}P_t(\mu_i(t))=(-1)^i\cdot2|a|,\quad i=1,...,p-1\ .\label{ecrit}\end{equation}
These conditions hold true for $t=0$ by construction, and we only need to verify that the Jacobian with respect to $\mu_1,...,\mu_{p-1}$ does not vanish. To this end, we observe that there exists a diffeomorphism
of a neighborhood of the point $(\mu_1^{(0)},...,\mu^{(0)}_{p-1})\in\R^{p-1}$ onto a neighborhood of the
point $(b^{(0)}_{p-2},...,b^{(0)}_0)\in\R^{p-1}$ sending the critical points of a polynomial
$\lambda^p+\widetilde b_{p-2}\lambda^{p-2}+...+\widetilde b_0$ to its coefficients. Then the Jacobian of
the left-hand side of the system (\ref{ecrit}) with respect to $\mu_1,...,\mu_{p-1}$ at $t=0$ turns to be
$$\det\left((\mu_i^{(0)})^j\frac{\partial b_j}{\partial \mu_i}\Big|_{t=0}\right)_{i=1,...,p-1}^{j=0,...,p-2}=
\pm\prod_{1\le i<j\le p-1}(\mu^{(0)}_i-\mu^{(0)}_j)\cdot\det\frac{D(\widetilde b_{p-2},...,\widetilde b_0)}
{D(\mu_1,...,\mu_{p-1})}\big|_{t=0}\ne0\ .$$
It follows from (\ref{ecrit}) that, for any $\theta\in\R/2\pi\Z$, the equation (\ref{ecrit1}) on $\sigma$
has $p$ real solutions (counting multiplicities) in the interval $|\sigma|<\lambda_0$, and we have exactly
$(p-1)(p+q)=\delta(C,z)-\ImBr(C,z)$ double roots as
$$\sigma=\mu_{2i-1}(t),\quad \cos((p+q)\theta-\theta_a)=-1,\quad 1\le i\le\frac{p}{2}\ ,$$
or $$\sigma=\mu_{2i}(t),\quad \cos((p+q)\theta-\theta_a)=1,\quad 1\le i\le\frac{p-1}{2}\ .$$
That is, family (\ref{ecrit2}) indeed describes a real morsification of $(C,z)$.

Note, that the real curve $\{F_t=0\}\subset\R B_{C,z}$ is an immersed circle lying in the
$\lambda_0t^{\frac{p+q}{p}}$-neighborhood of the ellipse $\rho=t$ and transversally intersecting in $2p$ points (counting multiplicities) with each real line through the origin.

\smallskip

\begin{figure}
\setlength{\unitlength}{0.8cm}
\begin{picture}(16,20)(-1,0)
\thinlines
\put(1,1){\vector(0,1){11.5}}\put(1,6.5){\vector(1,0){6.5}}
\put(9.5,1){\vector(0,1){11.5}}\put(9.5,6.5){\vector(1,0){6.5}}
\put(1,14.5){\vector(0,1){6}}\put(1,14.5){\vector(1,0){6.5}}
\put(9.5,14.5){\vector(0,1){6}}\put(9.5,14.5){\vector(1,0){6.5}}

\dashline{0.2}(1,11.5)(6,11.5)\dashline{0.2}(6,6.5)(6,11.5)
\dashline{0.2}(14.5,6.5)(14.5,11.5)\dashline{0.2}(3,14.5)(3,16.5)
\dashline{0.2}(1,16.5)(3,16.5)

\thicklines
\put(1,1.5){\line(0,1){5}}\put(1,1.5){\line(2,5){2}}
\put(1,6.5){\line(1,0){2}}\put(1,6.5){\line(1,1){5}}
\put(3,6.5){\line(3,5){3}}
\put(9.5,1.5){\line(0,1){10}}\put(9.5,1.5){\line(2,5){2}}
\put(11.5,6.5){\line(-2,5){2}}\put(9.5,11.5){\line(1,0){5}}
\put(11.5,6.5){\line(3,5){3}}
\put(3,16.5){\line(-2,3){2}}\put(3,16.5){\line(3,-2){3}}
\put(11.5,16.5){\line(-2,3){2}}\put(11.5,16.5){\line(3,-2){3}}
\put(9.5,14.5){\line(1,1){2}}\put(9.5,14.5){\line(1,0){5}}
\put(9.5,14.5){\line(0,1){5}}

\put(-0.4,19.4){$p+q$}\put(-0.4,11.4){$p+q$}
\put(8.1,19.4){$p+q$}\put(8.1,11.4){$p+q$}
\put(-0.7,1.5){$-p-q$}\put(7.8,1.5){$-p-q$}
\put(5.5,14){$p+q$}\put(14,14){$p+q$}
\put(5.5,6){$p+q$}\put(14,6){$p+q$}
\put(0.5,16.4){$p$}\put(2.9,14){$p$}
\put(11.5,6){$p$}\put(3,6){$p$}
\put(11.5,15){$T_1$}\put(10,16.5){$T_2$}
\put(1.5,5.2){$T'_1$}\put(2.3,7){$T'_2$}\put(10,7){$T$}
\put(7,14.7){$w$}\put(15.5,14.7){$w$}
\put(7,6.7){$u$}\put(1,20.7){$\widehat w$}\put(9.5,20.7){$\widehat w$}
\put(1,20.7){$\widehat w$}

\put(4,0.5){\rm (c)}\put(12.5,0.5){\rm (d)}
\put(4,13.3){\rm (a)}\put(12.5,13.3){\rm (b)}
\end{picture}
\caption{Patchworking a real morsification}\label{fcrit}
\end{figure}

{\bf(2)} Consider the general case.
By a coordinate change
$$(w,\widehat w)\mapsto
\left(w+\sum_{i\ge2}\alpha_i\widehat w^i,\ \widehat w+\sum_{i\ge2}\overline\alpha_i w^i\right)$$
one can bring $(C,z)$ to a strictly Newton non-degenerate form with the Newton diagram
$\Gamma(F)=[(p+q,0),p,p)]\cup[(p,p),(0,p+q)]$ in the
coordinates $w,\widehat w$ (see Figure \ref{fcrit}(a)), and with an equation
$$F(w,\widehat w)=(w\widehat w)^p-a\widehat w^{p+q}-\overline aw^{p+q}+\sum_{\renewcommand{\arraystretch}{0.6}
\begin{array}{c}
\scriptstyle{pi+qj>p(p+q)}\\
\scriptstyle{qi+pj>p(p+q)}\end{array}}a_{ij}w^i\widehat w^j=0\ ,$$
where $a\in\C^*$ and $a_{ij}=\overline a_{ji}$ for all $i,j$
(cf. (\ref{ecrit2})). We construct a real morsification of $(C,z)$ combining the result of the preceding step with the patchworking construction as developed in \cite[Section 2]{ST}.

Denote by $\Delta(F)$ the Newton polygon of $F(w,\widetilde w)$ and divide the domain under
$\Gamma(F)$ by the segment $[(0,0),(p,p)]$ into two triangles $T_1,T_2$
(see Figure \ref{fcrit}(b)). So, $\Delta(F)$, $T_1$, and $T_2$ form a convex subdivision of the
convex polygon
$\widetilde\Delta(F)=\conv(\Delta(F)\cup\{(0,0)\})$, i.e., there exists a convex piecewise linear
function $\nu:\widetilde\Delta(F)\to\R$ taking integral values at integral points and
whose linearity domains are $\Delta(F)$, $T_1$, and $T_2$. 
The overgraph
$\Graph^+(\nu)$ of $\nu$ is a three-dimensional convex lattice polytope, and we have a natural
morphism $\Tor(\Graph^+(\nu))\to\C$ whose fibers for $t\in\C^*$ are isomorphic to
$\Tor(\widetilde(F))$, and the central fiber is the union
$\Tor(\Delta(F))\cup\Tor(T_1)\cup\Tor(T_2)$. In the toric surface $\Tor(\Delta(F))$, we
have a curve $C=\{F(w,\widehat w)=0\}$, in the toric surfaces $\Tor(T_1)$ and $\Tor(T_2)$, we define curves
$$R_1=\{(w\widehat w-1)^p-\overline aw^{p+q}=0\}\quad\text{and}\quad
R_2=\{(w\widehat w-1)^p-a\widehat w=0\}\ ,$$
respectively. The complex conjugation interchanges the pairs $(\Tor(T_1),R_1)$ and
$(\Tor(T_2),R_2)$. Note that $R_1,R_2$ transversally intersect the toric divisors
$\Tor([(p,p),(p+q,0)])$ and $\Tor([(p,p),(0,p+q)])$ in the same points as
$C$. Furthermore, $R_1$, $R_2$ are rational curves intersecting the toric divisor
$\Tor([(0,0),(p,p)]=\Tor(T_1)\cap\Tor(T_2)$ in the same point $z_1$, where each of them has
a singular point of topological type $x^p+y^{p+q}=0$. To apply the patchworking
statement of \cite[Theorem 2.8]{ST}, we perform the weighted blow up
${\mathfrak X}\to\Tor(\Graph^+(\nu))$ of the point $z_1$ with the exceptional divisor
$E=\Tor(T)$, $T=\conv\{(p,0),(0,p+q),(0,-p-q)\}$ (see \cite[Figure 1]{ST}) being a part of the central fiber of
${\mathfrak X}\to\C$.

One can view this blow up via the refinement procedure developed in
\cite[Section 3.5]{Sh1}. Namely, we perform the toric coordinate change $u=w\widehat w$,
$v=w^{-1}$, transforming the triangles $T_1,T_2$ to $T'_1,T'_2$
as shown in Figure \ref{fcrit}(c), and respectively transforming the curves $R_1,R_2$ and the function $\nu$.
Note that this coordinate change defines an automorphism of the punctures real plane $\R^2\setminus\{0\}$.
Next we perform another coordinate change $u=u_1+1$, $v=v_1$, bringing the singular points of $R_1,R_2$ to the origin and transforming their Newton triangles $T'_1,T'_2$ into the edge $T''_1=[(p,0),(0,-p-q)]$ and the triangle
$T''_2=\conv\{(0,p+q),(p,0),(p+q,p+q)\}$, respectively (see Figure \ref{fcrit}(d)). The triangle
$T=\conv\{(0,-p-q),(0,p+q),(p,0)\}$ corresponds to the exceptional surface, in which we have to define a real curve by an equation with Newton triangle $T$, having the coefficients at the vertices determined by the equations of
$R_1$ and $R_2$ and having $(p-1)(p+q)=\delta(C,z)-\ImBr(C,z)$ real hyperbolic nodes. We just borrow the required curve from the special example studied in the first step. Namely, we do the above transformations with
the data given by (\ref{ecrit3}), and arrive at the curve given by a polynomial having
coefficient $a$ at $(0,p+q)$, coefficient $\overline a$ at $(0,-p-q)$, coefficient $1$ at $(p,0)$, and
coefficients $b_i^{(0)}$ at $(i,0)$, $i=0,...,p-2$.

To apply \cite[Theorem 2.8]{ST}, we have to verify the following transversality conditions:
\begin{itemize}\item for $i=1,2$, the germ at $R_i$ of the family of curves on the surface
$\Tor(T_i)$ in the tautological linear system that have a singularity of the topological type $x^p+y^{p+q}=0$ in a fixed position, is smooth of expected dimension;
\item the germ at $R$ of the family of curves on the surface $\Tor(T)$ in the tautological linear system that intersect the toric divisors $\Tor([(0,-p-q),(p,0)])$ and
    $\Tor([(p,0),(0,p+q)])$ in fixed points and have $(p-1)(p+q)$ nodes, is smooth of expected dimension.
\end{itemize} Both conditions are particular cases of the $S$-transvesality property, and they follow from the criterion in \cite[Theorem 4.1(1)]{ShT}. In the former case, one needs the inequality
$-R_iK_i>b$, where $K_i$ is the canonical divisor of the surface $\Tor(T_i)$, and $b$ a topological invariant of the singularity defined by
$$b(x^p+y^{p+q}=0)=\begin{cases}p+(p+q)-1,\quad & \text{if}\ q\not\equiv 1\mod p,\\
p+(p+q)-2,\quad &\text{if}\ q\equiv1\mod p\end{cases}$$
and the inequality holds, since $-R_iK_i=p+(p+q)+1$. In the latter case, one needs the inequality
$$R\cdot\Tor([(0,p+q),(0,-p-q)])>0$$ (nodes do not count in the criterion), which
evidently holds.

Thus, \cite[Theorem 2.8]{ST} yields the existence of an analytic equivariant deformation
of $F(w,\widehat w)$ defining in $\R B_{C,z}$ curves with $(p-1)(p+q)=\delta(C,z)-\ImBr(C,z)$ hyperbolic nodes.
\proofend

\begin{lemma}\label{lNDD1}
Let a real singularity $(C,z)$ with exactly two tangent lines $L,\overline L$ be admissible along its tangent lines. Then $(C,z)$ possesses a real morsification.
\end{lemma}

{\bf Proof.}
We apply construction presented in the proof of Lemmas \ref{lsmooth} and \ref{lNND} for the bunch of smooth branches a nd for pairs of singular complex conjugate branches separately, and we shall show that, for
any two pairs $(C_1,\overline C_1)$, $(C_2,\overline C_2)$ of complex conjugate branches of $(C,z)$, their
divides intersect in $2(C_1\cdot C_2)+2\mt C_1\cdot\mt C_2$ (real) points.

For $C_1,C_2$ smooth this follows from Lemma \ref{lsmooth}. In other situations, we can assume that $C_1\cup C_2$
(and $\overline C_1\cup\overline C_2$) form a strictly Newton non-degenerate singularity so that $C_1$ os of topological type
$x^p+y^q=0$ with $2\le p<q$, $(p,q)=1$, and $C_2$ is of topological type $x^{p'}+y^{q'}=0$ with
$1\le p'<q'$, $(p',q')=1$.

If $q/p=q'/p'$, then $p=p'$, $q=q'$, and hence $C_1\cup\overline C_1$ and $C_2\cup C_2$ are given by
$$F(w,\widehat w)=(w\widehat w)^p-a\widehat w^{p+q}-\overline aw^{p+q}+\sum_{\renewcommand{\arraystretch}{0.6}
\begin{array}{c}
\scriptstyle{pi+qj>p(p+q)}\\
\scriptstyle{qi+pj>p(p+q)}\end{array}}a_{ij}w^i\widehat w^j=0\ ,$$
and
$$F'(w,\widehat w)=(w\widehat w)^p-a'\widehat w^{p+q}-\overline a'w^{p+q}+\sum_{\renewcommand{\arraystretch}{0.6}
\begin{array}{c}
\scriptstyle{pi+qj>p(p+q)}\\
\scriptstyle{qi+pj>p(p+q)}\end{array}}a'_{ij}w^i\widehat w^j=0\ ,$$
respectively, where $a,a',a-a'\in\C^*$. The patchworking construction in the second step of
the proof of Lemma \ref{lNND} can be applied to both the pairs of the branches simultaneously, and the considered question on the intersection of the divides reduces then to the intersection of the curves $R,R'$ in the toric surface $\Tor(T)$, $T=\conv\{(0,-p-q),(p,0),(0,p+q)\}$. The real parts $\R R,\R R'$ of
these curves, in suitable coordinates $\sigma>0$, $\theta\in\R/2\pi\Z$ are given by
$$\sigma^p+\sum_{i=p-2}^0b_i^{(0)}\sigma^i=2|a|\cos((p+q)\theta-\theta_a),\quad
\sigma^p+\sum_{i=p-2}^0b_i^{(0)}\sigma^i=2|a'|\cos((p+q)\theta-\theta_{a'})\ ,$$
respectively. The number of their (real) intersection points is $p$ times the number of
solutions of the equation
$$|a|\cos((p+q)\theta-\theta_a)=|a'|\cos((p+q)\theta-\theta_{a'}),\quad \theta \in\R/2\pi\Z\ .$$
The latter number is $2(p+q)$, and hence the total number of intersection points is
$$2p(p+q)=2pq+2p^2=2(C_1\cdot C_2)+2\mt C_1\cdot\mt C_2$$ as required.

Suppose that $\tau=\frac{q'}{p'}-\frac{q}{p}>0$. Then $C_1\cup\overline C_1$ and $C_2\cup C_2$ are given by
$$F(w,\widehat w)=(w\widehat w)^p-a\widehat w^{p+q}-\overline aw^{p+q}+\sum_{\renewcommand{\arraystretch}{0.6}
\begin{array}{c}
\scriptstyle{pi+qj>p(p+q)}\\
\scriptstyle{qi+pj>p(p+q)}\end{array}}a_{ij}w^i\widehat w^j=0\ ,$$
and
$$F'(w,\widehat w)=(w\widehat w)^{p'}-a'\widehat w^{p'+q'}-\overline a'w^{p'+q'}+\sum_{\renewcommand{\arraystretch}{0.6}
\begin{array}{c}
\scriptstyle{p'i+q'j>p'(p'+q')}\\
\scriptstyle{q'i+p'j>p'(p'+q')}\end{array}}a'_{ij}w^i\widehat w^j=0\ ,$$
respectively. Along the construction of Lemmas \ref{lsmooth} and \ref{lNND}, we substitute in the above equations
$(w\widehat w-t^2)^p$ for $(w\widehat w)^p$ and $(w\widehat w-t^2)^{p'}$ for
$w\widehat w)^{p'}$, respectively, then make the same rescaling $(w,\widehat w)\mapsto(tw,t\widehat w)$. Next, we
pass to the real coordinates $\sigma,\theta$ via
$$\rho^2=w\widehat w=1+t^{\frac{q-p}{p}}\sigma, \quad w=\rho\exp(\sqrt{-1}\theta),\ \widehat w=\rho\exp(-\sqrt{-1}\theta)\ ,$$
(adapted to the pair $p,q$, not $p',q'$ !). Then the real morsification of $C_1\cup\overline C_1$ is given by
$$\sigma^p+\sum_{i=p-2}^0b_i^{(0)}\sigma^i=2|a|\cos((p+q)\theta-\theta_a)+O(t^{\frac{1}{p}})\ ,$$
while the real morsification of $C_2\cup\overline C_2$ is given by
$\sigma^{p'}=O(t^{\\tau})$. The divide of the real morsification of $C_2\cup\overline C_2$ is the circle immersed into
the $O(t^{\frac{1}{p'}})$-neighborhood of the level line $\sigma=0$ in the annulus
$\{(\sigma,\theta)\in(-\lambda_0,\lambda_0)\times(\R/2\pi\Z)\}$ so that the normal projection onto the circle $\sigma=0$
is a $p'$-fold covering. Hence, this divide intersects with the divide of the real morsification of $C_1\cup\overline C_1$ in
$2p'(p+q)=2p'q+2p'p=2(C_1\cdot C_2)+2\mt C_1\cdot\mt C_2$ real points.

The case of $\tau=\frac{q}{p}-\frac{q'}{p'}<0$ can be considered in the same manner.
\proofend

\subsubsection{The case of several pairs of complex conjugate tangents}\label{sec-sev}
Suppose now that $(C,z)$ has $r\ge2$ pairs of complex conjugate tangent lines
$$L_i=\{x+(\alpha_i+\beta_i\sqrt{-1})y=0\},
\quad\overline L_i=\{x+(\alpha_i-\beta_i\sqrt{-1})y=0\},\quad i=1,...,r\ ,$$
where $\alpha_i,\beta_i\in\R$, $\beta_i\ne0$ for all $i=1,...,r$.
Set $$w_i=x+(\alpha_i+\beta_i\sqrt{-1})y,\quad\widehat w_i=x+(\alpha_i-\beta_i\sqrt{-1})y,
\quad i=1,,,.,r\ .$$ Equations $\rho_i^2:=w_i\widehat w_i=\const>0$, $i=1,...,r$, define distinct ellipses in $\R^2$,
and there are $\gamma_1,...,\gamma_r>0$ such that each two ellipses $\rho_i^2=\gamma_i$,
$\rho_j^2=\gamma_j$, $1\le i<j\le r$, intersect in four (real) points,
and all $2r(r-1)$ intersection points are distinct.

For any $i=1,...,r$, we introduce a real singularity $(C^{(i)},z)$ formed by the union of all the branches of $(C,z)$ tangent either to $L_i$, or to $\overline L_i$, and then construct a real morsification of $(C^{(i)},z)$ following the procedure of
Section \ref{sec-smooth}, in which $t$ should be replaced by $t\sqrt{\gamma_i}$. For a given $t>0$, the divide of
this morsification lies in $O(t^{>2})$-neighborhood of the ellipse $\rho_i^2=\gamma_it^2$, and it is the union of several
immersed circles so that the normal projection onto the ellipse is a covering of multiplicity
$\frac{1}{2}\mt(C^{(i)},z)$. Hence,
the divides of the morsifications of $(C^{(i)},z)$ and $C^{(j)},z)$, $1\le i<j\le r$, intersect in
$\mt C^{(i)}\cdot\mt C^{(j)}$ real points. So, in total the union of all $r$ divides contains
$$\sum_{i=1}^r\left(\delta(C^{(i)},z)-\ImBr(C^{(i)},z)\right)+\sum_{1\le i<j\le r}
(C^{(i)}\cdot C^{(j)})_z=\delta(C,z)-\ImBr(C,z)$$
real hyperbolic nodes.

\subsection{Proof of Theorem \ref{t1}: general case}\label{sec-theorem1}
Suppose now that $(C,z)$ is a real singularity satisfying hypotheses of Theorem \ref{t1}.
Denote by $(C^{re},z)$, resp. $(C^{im},z)$, the union of the branches of $(C,z)$ that have real, resp.
complex conjugate tangents.

If $C^{re}=\emptyset$, the existence of a real morsification follows from the results of
Sections \ref{sec-smooth} and \ref{sec-sev}.
Assume that $C^{re}\ne\emptyset$, and it contains only smooth branches. We settle this case by induction on
$\Delta_3(C,z)$, the maximal $\delta$-invariant of a subgerm of $(C^{re},z)$ having a unique tangent line.
If $\Delta_3(C,z)=0$, then all branches
of $(C^{re},z)$ are smooth real and transversal to each other. Then we first construct a real morsification of $(C^{im},z)$ as in
Sections \ref{sec-smooth} and \ref{sec-sev} with $t>0$ chosen so small that each branch of $(C^{re},z)$
intersects the divide of the morsification of $(C^{im},z)$ in $\mt(C^{im},z)$ real points.
Then we
slightly shift the branches of $(C^{re},z)$ in general position keeping the above
real intersection points and obtaining additional $\delta(C^{re},z)$ hyperbolic nodes
as required. In the case of $\Delta_3(C,z)>0$, we blow up the point $z$ and consider the strict
transform of $(C^{re},z)$, which consists of germ $(C_i,z_i)$ with real centers $z_i$ on the exceptional
divisor $E$. Clearly, for each germ $(C_i\cup E,z_i)$, its branches with real tangents are smooth
and transversal to $E$,
and, furthermore, $\Delta_3(C_i\cup E,z_i)<\Delta_3(C,z)$ for all $i$. Hence, there are real morsifications of the germs $(C_i\cup E,z_i)$, in which we cam assume the germs $(E,z_i)$ to be fixed. Then we blow down
$E$ and obtain a deformation of $(C^{re},z)$ with $\mt(C^{re},z)$ real smooth transversal branches at $z$
and additional $\delta(C^{re},z)-\ImBr(C^{re},z)-\frac{1}{2}\mt(C^{re},z)(\mt(C^{re},z)-1)$
real hyperbolic nodes (cf. computations in Section \ref{total}(1)). Returning back the subgerm $(C^{im},z)$,
we obtain a real singularity at $z$ with $\Delta_3=0$, and thus, complete the construction of a real morsification of $(C,z)$ as in the beginning of this paragraph.

Now we get rid of all extra restrictions on
$(C^{re},z)$ and prove the existence of a real morsification of $(C,z)$ by
induction on $\Delta_4(C,z)$, which is the $\delta$-invariant of the union of singular branches
of $(C^{re},z)$. The preceding consideration serves as the base of induction. The induction step is precisely the same, and we only notice that (in the above notations)
$\max\Delta_4(C_i\cup E,z_i)<\Delta_4(C,z)$.

The proof of Theorem \ref{t1} is completed.

\section{Real morsifications and Milnor fibers}\label{sec-reg}

\subsection{A'Campo surface and Milnor fiber}\label{sec-mil}
In \cite[Section 3]{AC1}, A'Campo constructs the link of a divide of a real morsification of a singularity
(which we call {\bf A'Campo link}).
This link is embedded into the $3$-sphere, the boundary of the Milnor ball, and the fundamental result by
A'Campo \cite[Theorem 2]{AC1} states that it is isotopic to the link of the given singularity in the $3$-sphere.
In this section, we discuss a somewhat stronger isotopy. Namely, in \cite[Section 3]{AC1}, A'Campo
associates with a real morsification a surface (which we call {\bf A'Campo surface}), whose boundary is the A'Campo link. It is natural to ask whether
the pair (A'Campo surface, A'Campo link) is isotopic to the pair (Milnor fiber, its boundary).

In \cite[Page 22]{AC1}, A'Campo conjectures a certain transversality condition for the known morsifications
that ensure the discussed transversality. 
Here we prove this transversality condition for all morsifications constructed in Section \ref{sec-exi}.
We also show that the spoken transversality condition may fail even for morsifications of simple singularities.
Hence, the question on the isotopy between the A'Campo surface and the Milnor fiber remains open in a general case.

Let $(C,0)\subset\C^2$ be a real singularity given by an equivariant analytic equation $f(x,y)=0$.
Following \cite[Section 3]{AC1}, we replace the standard Milnor ball $B(C,0)$ by the bi-disc
$B(0,\rho_0):=\{u+v\sqrt{-1}\in\C^2\ :\ u,v\in D(0,\rho_0)\subset\R^2\}$, where $\rho_0>0$ and $\C^2=\R^2\oplus\R^2\sqrt{-1}$.
It is easy to verify that
$\partial B(0,\rho)$ transversally intersects with $C$ for each $0<\rho\le\rho_0$ if $\rho_0$ is small enough,
and we assume this further on. For $\xi\in\C$ with $0<|\xi|\ll1$ all curves
$M_\xi=\{f(x,y)=\xi\}\subset B(0,\rho_0)$ are smooth and transversally intersect $\partial B(0,\rho_0)$. They are called {\bf Milnor fibers} of the given singularity
$(C,0)$. Respectively, the links $LM_\xi=M_\xi\cap\partial B(0,\rho_0)$ are isotopic
in the sphere $\partial B(0,\rho_0)$ to the link $L(C,z)=C\cap\partial
B(0,\rho_0)$ of the singularity $(C,z)$, and the pairs $(M_\xi,LM_\xi)$, $0<|\xi|\ll1$, are isotopic in
$(B(0,\rho_0),\partial B(0,\rho_0))$.

Introduce the family of bi-discs
$$B'_\rho(0,\rho_0)=\{u+v\sqrt{-1}\in\C^2\ :\ u\in D(0,\rho_0),\ v\in D(0,\rho)\},\quad 0<\rho\le\rho_0\ .$$ By definition,
$B'_{\rho_0}(0,\rho_0)=B(0,\rho_0)$. Let $C_t=\{f_t(x,y)=0\}$, $0\le t\le t_0$, $f_0=f$, be a real morsification of $(C,0)$
defined in $B(0,\rho_0)$. Without loss of generality, we can assume that $C_t$ intersects with $\partial B(0,\rho_0)$ transversally for all $0\le
t\le t_0$.

We have two families of singular surfaces in $B(0,\rho_0)$:
\begin{itemize}\item $F(\rho)=C_{t_0}\cap B'_\rho(0,\rho_0)$, $0\le\rho\le\rho_0$,
\item $R(\rho)=\{u+v\sqrt{-1}\in B'_\rho(0,\rho_0)\ :\ u\in\R C_{t_0},\ v\in T_u\R C_{t_0},\ v\in D(0,\rho)\}$, $0\le\rho\le\rho_0$
(here $\R C_{t_0}\subset D(0,\rho_0)$ is an immersed real analytic curve with nodes, and at each node $u\in\R C_{t_0}$
we understand $T_u\R C_{t_0}$ as the union of the tangent lines to the branches centered at $u$).
\end{itemize} Denote $LF(\rho)=F(\rho)\cap\partial B'_\rho(0,\rho_0)$ and
$LR(\rho)=R(\rho)\cap\partial B'_\rho(0,\rho_0)$ for all $0<\rho\le\rho_0$.

\begin{lemma}\label{l15}[cf. \cite{AC1}, Theorem 2]
(1) The set $LR(\rho)$ is a link in the sphere $\partial B'(\rho)$ for any $0<\rho\le\rho_0$. The set
$LF(\rho)$ is a link in the sphere $\partial B'_\rho(0,\rho_0)$ for all but finitely many values $\rho\in(0,\rho_0]$.
Furthermore, $LF(\rho_0)$ is a link equivariantly isotopic in $\partial B(0,\rho_0)$ to the singularity link
$L(C,z)$.

(2) There exists $\rho'=\rho'(t_0)$ such that the links $LF(\rho')$ and
$LR(\rho')$ are equivariantly isotopic in $\partial B'_{\rho'}(0,\rho_0)$, and the pairs
$(F(\rho'),LF(\rho'))$ and $(R(\rho'),LR(\rho'))$ are equivariantly isotopic in
$(B'_{\rho'}(0,\rho_0),\partial B'_{\rho'}(0,\rho_0))$.
\end{lemma}

{\bf Proof.} The first statement is straightforward. The second one immediately follows from the fact that
$F(\rho)$ and $R(\rho)$ are immersed surfaces having the same real point set with the same tangent planes along it.
\proofend

For $\eta>0$ small enough, the algebraic curves $$F^{sm}(\rho)=\{f_{t_0}(x,y)=\eta\}\cap B'_\rho(0,\rho_0)$$ are smooth for all $\rho'(t_0)\le\rho\le\rho_0$,
and each of them is obtained from $F(\rho)$ by a small deformation in a neighborhood $U_u$ of each node $u\in \R C_{t_0}$ that
replaces two trasversally intersecting discs with a cylinder. Respectively, for all $\rho'(t_0)\le\rho\le\rho_0$,
we define $C^\infty$-smooth equivariant {\bf A'Campo surfaces}
$R^{sm}(\rho)\subset B'_\rho(0,\rho_0)$, obtained from $R(\rho)$ by replacing $R(\rho)\cap U_u$ with the cylinder $F^{sm}(\rho)\cap U_u$
smoothly attached to $R(\rho)\setminus U_u$ for each node $u\in\R C_{t_0}$.

If $\xi\in\C\setminus\{0\}$ with $|\xi|$ small enough, then the intersections
$M_\xi\cap\partial B'(\rho)$ are transversal for all $\rho'(t_0)\le\rho\le\rho_0$. We would like to address

\medskip

{\bf Question.} {\it Is the pair $(R^{sm}(\rho_0),LR(\rho_0))$ isotopic to $(M_\xi,LM_\xi)$ in $(B(0,\rho_0),\partial
B(0,\rho_0))$, or, equivalently, is the pair
$(R^{sm}(\rho'(t_0)),LR(\rho'(t_0)))$ isotopic to $(M_\xi\cap B'(\rho'(t_0)),M_\xi\cap\partial B'_{\rho'(t_0)}
(0,\rho_0)$ in
$(B'_{\rho'(t_0)}(0,\rho_0),\partial B'_{\rho'(t_0)}(0,\rho_0))$?}

\medskip

This seems to be stronger that Lemma \ref{l15}. We would like to comment on this question more.
Since $(F^{sm}(\rho_0),F^{sm}(\rho_0)\cap\partial B(0,\rho_0))$ is isotopic to $(M_\xi,LM_\xi)$ in
$(B(0,\rho_0),\partial B(0,\rho_0))$, and, by Lemma \ref{l15}, $(F^{sm}(\rho'(t_0)),F^{sm}(\rho'(t_0))\cap\partial
B'_{\rho'(t_0)}(0,\rho_0))$ is (equivariantly) isotopic to $(R^{sm}(\rho'(t_0)),LR(\rho'(t_0)))$ in
$B'_{\rho'(t_0)}(0,\rho_0),\partial B'_{\rho'(t_0)}
(0,\rho_0))$, the answer to the above Question would be {\bf yes}, if we could prove one of the following
claims. Observe that the closure of $R^{sm}(\rho_0)\setminus R^{sm}(\rho'(t_0))$ as well as the closure of $F^{sm}(\rho_0)\setminus
F^{sm}(\rho'(t_0))$ is the disjoint union of pairs of discs (corresponding to real branches of $(C,z)$) and cylinders
(corresponding to pairs of complex conjugate branches of $(C,z)$), and the former surface defines a cobordism of
$LR(\rho_0)$ and $LR(\rho'(t_0))$ trivially fibred over $[\rho'(t_0),\rho_0]$. So the requested claims are

\begin{enumerate}\item[(A)] The surface $Closure(F^{sm}(\rho_0)\setminus
F^{sm}(\rho'(t_0)))$ defines a trivial cobordism of $F^{sm}(\rho_0)\cap\partial B(0,\rho_0)$ and
$F^{sm}(\rho'(t_0))\cap\partial B'_{\rho'(t_0)}(0,\rho_0)$.
\item[(B)] The intersections $C_t\cap\partial B'_{\rho'(t_0)}(0,\rho_0)$ are trasversal for all $0\le t\le t_0$.
\end{enumerate}

Claim (A) seems to be open in general so far, and it is proved in
\cite{P} for morsifications of totally real singularities obtained by the blowing up construction
as in \cite{AC} (see also \cite[Theorem 5.2]{CP}).
Claim (B) is formulated in \cite[Page 22]{AC1} as a conjecture again for the morsifications
of totally real singularities constructed in \cite{AC}. However, in general, it does not hold:

\begin{proposition}\label{l16}
The totally real singularity $(C,z)$ given by $y^2-x^{2n}=0$, $n\ge4$, possesses a real morsification $C_t$, $0\le t\le t_0$ such that
for arbitrary $0<\rho<\rho_0$ and $0<t<t_0$, there exist $0<\rho'<\rho$ and $0<t'<t$ for which the intersection
of $C_{t'}$ and $\partial B'_{\rho'}(0,\rho_0)$ is not transversal.
\end{proposition}

{\bf Proof.}
We have $\partial B'_\rho(0,\rho_0)=\left(\partial D(0,\rho_0)\times D(0,\rho)\right)\cup \left(D(0,\rho_0)\times\partial D(0,\rho)\right)$. The intersection of $C_t$ with
$\partial D(0,\rho_0)\times D(0,\rho)$ is transversal for any real morsification of $(C,z)$. On the other hand, the intersection
of $C_t$ with $D(0,\rho_0)\times\partial D(0,\rho)$ is not transversal at some point $p=u+v\sqrt{-1}\in D(0,\rho_0)\times\partial D(0,\rho)$
if and only if the tangent line to $C_t$ at this point has a real slope. Indeed, if $C_t$ is given in a neighborhood of $p$ by
$y=\varphi(x)$, then the lack of transversality of the intersection of $C_t$ and $D(0,\rho_0)\times\partial D(0,\rho)$ at $p$ can be expressed as
$$\Ima\frac{d\varphi}{dx}\big|_p\cdot v_2=v_1-\Rea\frac{d\varphi}{dx}\big|_p\cdot v_2=0,\quad\text{where}\ v=(v_1,v_2)\ne0\ ,$$
and hence $\Ima\frac{d\varphi}{dx}\big|_p=0$. In other words, the lack of transversality means the existence of a real slope tangent line
to $C_t$ at a non-real point.

Now we define
$$C_t=\left\{(y-tx^2)^2-\prod_{k=1}^n(x-kt)^2=0\right\},\quad 0\le t\le t_0,\quad 0<t_0\ll1\ .$$
The real point set of $C_t$ consist of two branches $y=tx^2\pm\prod_{k=1}^n(x-kt)$ transversally intersecting in $n$ points, and
hence it is a real morsification. It is easy to compute that the branch $y=tx^2+\prod_{k=1}^n(x-kt)$ has $n-2$ tangent lines with
the zero slope at the points
$$x_i(t)=\lambda^i\left(\frac{2}{n}\right)^{1/(n-2)}t^{1/(n-2)}(1+O(t^{>0})),\quad i=0,...,n-3\ ,$$ where $\lambda^{n-2}=-1$ is a primitive root of unity. Thus, we obtain at least $n-3$ zero slope tangents at imaginary points. Since $x_i(t)\to0$ as $t\to0$,
the statement of Proposition follows.
\proofend

\subsection{Real Milnor morsifications}
We say that a real morsification of a real singularity $(C,z)$ is a {\bf real Milnor morsification} if
in the notation of Section \ref{sec-mil}, the pair $(R^{sm}(\rho_0),LR(\rho_0))$ is
isotopic to $(M_\xi,LM_\xi)$ in $(B'_\rho(z,\rho_0),\partial
B'_\rho(z,\rho_0))$ for some $0<\rho\le\rho_0$.

\begin{theorem}\label{ap2}
Any isolated real plane curve singularity satisfying the hypotheses of Theorem \ref{t1}
admits a real Milnor morsification.
\end{theorem}

{\bf Proof.} We prove the theorem by establishing Claim (B) formulated in the preceding section.

Let $(C,z)$ be a real singularity as in Theorem \ref{t1}. Applying a suitable local diffeomorphism, we
can assume that $(C,z)$ does not contain (segments of) straight lines, and
hence $(L\cdot C)_z<\infty$ for any line $L$ through $z$. Denote by $\Lambda$ the union of
all real tangent lines to $(C,z)$ at $z$.
Under the assumption made, we apply the construction used in the
proof of Theorem \ref{t1} and obtain a real morsification of
$(C\cup\Lambda,z)$, in which $\Lambda$ remains fixed. Then we get rid of
$\Lambda$ and obtain a real morsification $C_t$, $0\le t\le t_0$, of $(C,z)$. We shall show that
it is a real Milnor morsification (possibly replacing $t_0$ with a smaller positive number).

As noticed in the proof of Proposition \ref{l16}, the required property
is equivalent to the absence of non-real lines with real slopes tangent to $C_t$, $0\le t\le t_0$.

Our first observation is

\begin{lemma}\label{l7}
Let $(C,z)$ be a real singularity, $L$ a real line passing through $z$
and intersecting $(C,z)$ only at $z$ (in the Milnor ball), with a finite multiplicity
$(L\cdot C)_0$. Denote by
${\mathcal P}_L$ the germ of the pencil of the lines parallel to $L$ and by
$\R{\mathcal P}_L$ its real point set. Let $C_t$, $0\le t<\eps$, be
a real morsification of $(C,z)$ as above, and let $C_t$ and $L_t$ intersect in
$(L\cdot C)_z$ real points
for any $t\in(0,\eps)$. Then each line $L'\in{\mathcal P}_L\setminus\R{\mathcal P}_L$
intersects each element $C_t$, $0< t<\eps$, transversally.
\end{lemma}

{\bf Proof.} Let $C'$ be a Milnor fiber. Then the lines of
${\mathcal P}_L$ in total are tangent to $C'$ in $\kappa(C,z)+(L\cdot C)_z-\mt(C,z)$ points,
where $\kappa(C,z)$ is the class of the singularity $(C,z)$ (see, for example,
\cite[Section I.3.4]{GLS} for details). Since, for a node, $\kappa=2$, and in general
$\kappa(C,z)=2\delta(C,z)+\mt(C,z)-\Br(C,z)$, we get that
the lines of ${\mathcal P}_L$ in total are tangent to $C_t$ in
$$\kappa(C,z)+(L\cdot C)_z-\mt(C,z)-2(\delta(C,z)-\ImBr(C,z))=
(L\cdot C)_z-\ReBr(C,z)$$ points. It follows that
\begin{itemize}\item $L$ intersects the morsification $C_{i,t}$ of
any real branch $(C_i,z)$ of $(C,z)$ in $(L\cdot C_i)_z$
real points, while the real point set $\R C_{i,t}$
of $C_{i,t}$ is an immersed segment; that is, $L$ cuts $\R C_{i,t}$ into
$(L\cdot C_i)+1$ immersed segments, among all but two have both endpoints on $\R L$;
hence, varying $L$ in $\R{\mathcal P}_L$, we encounter
at least $(L\cdot C_i)_z-1$ real tangency points;
\item $L$ intersects the morsification $C_{j,t}$ of a pair of complex conjugate branches
$(C_j,z)$, $\overline C_j,z)$ of $(C,z)$ in $2(L\cdot C_i)_z$ real points, and hence it cuts
$\R C_{j,t}$ (which is an immersed circle) into $2(L\cdot C_i)_z$ immersed segments, whose all endpoints lie
on $\R L$, and hence, varying $L$ in $\R{\mathcal P}_L$, we encounter at least
$2(L\cdot C_i)_z$ real tangency points.
\end{itemize} The claim of Lemma follows.
\proofend

Remark that, under conditions of Lemma \ref{l7}, there is an open neighborhood $U_L$ of $L$ in the dual plane
$\PP^{2,\vee}$ such that all non-real curves with real slopes intersect each curve $C_t$, $0<t<\eps$, transversally.
Thus, Theorem \ref{ap2} follows from

\begin{lemma}\label{al7}
For any real line $L$ through $z$, there exist $0<\rho\le\rho_0$ satisfying the
following conditions
\begin{itemize}\item $L\cap C\cap B'_\rho(z,\rho_0)=\{z\}$; \item for some $\eps>0$, $L$ intersects with
any curve $C_t$, $0<t<\eps$, in $(L\cdot C)_z$ real points (counting multiplicities).
\end{itemize}
\end{lemma}

{\bf Proof.} Let $L_1,...,L_k$ be all real tangent lines to $(C,z)$ at $z$.
Write $(C,z)=\bigcup_i(C_i,z)$, where $(C_i,z)$ either has a unique (real) tangent line, or a pair of complex conjugate tangent lines, and $(C_i,z)$, $(C_j,z)$ have no tangent in common as $i\ne j$. We can consider
morsifications of $(C_i,z)$ separately.

Suppose that $(C_i,z)$ has a pair of complex conjugate tangent lines. The morsification of $(C_i,z)$ constructed in Section \ref{sec-smooth} is such that the
real point set of $C_t$, $0<t<\eps$, consists of one or several immersed circles going
in total $\frac{1}{2}\mt(C_i,z)$ times around $z$, and hence $L$ (which is transversal to $(C_i,z)$, i.e.
$(L\cdot C_i)_z=\mt(C_i,z)$) intersects
any curve $C_t$ in $\mt(C_i,z)$ real points (counting multiplicities).

Suppose that $(C_i,z)$ has a unique (real) tangent line $L_z$, and $L\ne L_z$. Then $(L\cdot C_i)_z=\mt(C_i,z)$.
The smooth real branches of $(C_i,z)$ are deformed in any morsification so that they remain transversal to $L$ and intersect $L$ at one real point. For $(C'_i,z)$, the union of the other branches of $(C_i,z)$, the construction of a morsification presented in Section \ref{sec-theorem1} goes inductively.
Namely, we blow up $z$, construct a morsification of the
strict transform of $(C_i,z)$ united with the exceptional divisor and then blow down the exceptional divisor.
Elements of this intermediate deformation have $\mt(C'_i,z)$ smooth real branches centered at $z$,
all transversal to $L$, and in any further deformation they intersect with $L$ in $\mt(C'_i,z)$ real points.

If $(C_i,z)$ has a unique (real) tangent line $L_z$, and $L=L_z$, the statement follows from the construction.
\proofend

\section{A'Campo-Gusein-Zade diagrams and topology of singularities}\label{sec3}

\subsection{$\AG$-diagrams of real morsifications}\label{sec-ag}
L. Balke and R. Kaenders proved \cite[Theorem 2.5 and Corollary 2.6]{BK} that
the A'Campo-Gusein-Zade diagram (briefly, $\AG$-diagram) associated with a morsification of a totally real singularity determines the complex topological type of the given singularity. Here we extend this result to real morsifications of arbitrary real singularities. We
get rid of the requirement for morsifications to define a partition (see Section \ref{sec1}
and \cite[Definition 1.2]{BK}) and prove that an $\AG$-diagram determines the topological type
of the singularity as well as some additional information on its real structure.

Let us recall definitions from \cite{AC3} and \cite{BK}.

A subset $D$ of a closed disc $\DD\subset \mathbb{R}^2$ is called a {\bf connected divide} if it is the image of an immersion of a disjoint union $\Sigma\ne\emptyset$ of a finite number of segments $I=[0,1]$ and circles $S^1$ satisfying the following conditions:
	\begin{itemize}
		\item the set of the endpoints of all the segments in $\Sigma$ is injectively mapped to $\partial\DD$,
whereas the other points of $\Sigma$ are mapped to the interior of $\DD$;
		\item each point of the complement $D\setminus\Sing(D)$ to a finite set $\Sing(D)$ has a
unique preimage in $\Sigma$, each point of $\Sing(D)$ is a transversal intersection of
two smooth local branches;
\item the images of any two connected components of $\Sigma$ intersect each other.
	\end{itemize}

Note that $\Sigma$ is uniquely determined by $D$. The image of any connected component of $\Sigma$
is a divide, which is called a {\bf branch} of the divide $D$.

The divide of a real morsification of a real singularity placed in the real Milnor disc (see Section \ref{sec1}) is a connected divide
in the above sense.

Connected components of $\overline D\setminus D$ and of $D\setminus\Sing(D)$, disjoint from $\partial\DD$, are called inner components. Clearly, each inner component of
$\DD\setminus D$ is homeomorphic to an open disc, and each inner component
of $D\setminus\Sing(D)$ is homeomorphic either to an open interval, or to $S^1$ if $D\simeq S^1$.

It is straightforward that the set
$\pi_0(\DD\setminus D)$ of the connected components of $\DD\setminus D$ can be $2$-colored, i.e.,
there exists a function $\pi_0(\DD\setminus D)\to\{\pm1\}$ such that the components, whose boundaries
intersect along one-dimensional pieces of $D$, have different signs, and there are precisely two functions like that (cf. \cite[Proposition 1.4]{BK}). Fix a $2$-coloring
$s:\pi_0(\DD\setminus D)\to\{\pm1\}$. The {\bf A'Campo-Gusein-Zade diagram} ($\AG${\bf -diagram}) of a connected divide $D$ is a $3$-colored graph $\AG(D)=(V,E,c)$ such that
\begin{itemize}\item the set $V$ of its vertices is in one-to-one correspondence with
the disjoint union of $\Sing(D)$ (the set of $\bullet$-vertices in the notation
of \cite{BK}) and the set $\pi^{inn}_0(\DD\setminus D)$ of the
inner components of $\DD\setminus D$ (the $\oplus$-vertices and $\ominus$-vertices
in the notation of \cite{BK} in accordance with the chosen coloring);
\item two distinct vertices $K_1,K_2\in\pi^{inn}_0(\DD
\setminus D)$ such that $\partial K_1\cap\partial K_2\setminus\Sing(D)\ne\emptyset$ are joined by $k$ edges, where $k$ is the number of inner components of
$D\setminus\Sing(D)$ inside $\partial K_1\cap\partial K_2$;
\item two vertices $K\in\pi^{inn}_0(\DD\setminus D)$ and $p\in\Sing(D)$ such that
$p\in\partial K$ are joined by $k$ edges, where $k$ is the number of components of the intersection
of $K$ with a small disc centered at $p$ (clearly, here $k=1$ or $2$);
\item the $3$-coloring $c:V\to\{\pm1,0\}$ is defined by
$c(K)=s(K)$, $K\in\pi^{int}_0(\DD\setminus D)$, and $c(p)=0$, $p\in\Sing(D)$.\end{itemize}
Comparing with \cite[Definition 1.5]{BK}, we admit multi-graphs, i.e., vertices can be joined by several edges, while this is excluded in \cite[Definition 1.5]{BK} by the partition requirement.
On the other hand, there are no loops. By construction, the $\AG$-diagram can be embedded into
$D$ (cf. \cite[Remark in page 43]{BK}).

The $\AG$-diagram associated with the divide of a real morsification of a real singularity is simply called an $\AG$-diagram of that singularity.

\subsection{$\AG$-diagram determines the weak real topological type of a singularity}
The topological type
of a real singularity $(C,z)$ is its equivalence class up to a homomorphism of the Milnor ball, and it is known \cite{Bra,Z} (see
also \cite[Section 8.4]{BrK}) that the topological type of a given singularity is determined by the
collections of Puiseux pairs of its branches and by pairwise intersection numbers
of the branches. We introduce the {\bf weak real topological type} of $(C,z)$ to be the topological type enriched with the following information:
\begin{itemize}\item indication of real branches and pairs of complex conjugate branches;
\item the cyclic order of real branches, that is, if $(C,z)$ has $k\ge1$ real branches, we number them somehow
and introduce the cyclic order on the
multiset $\{1,1,2,2,...,k,k\}$ induced by the position of the $2k$ intersection points
of the real branches with the circle $\partial\R B_{C,z}$ and defined up to reversing the orientation of $\partial\R B_{C,z}$
and renumbering the topological types of the real branches, their mutual intersection multiplicities and their intersection multiplicities with
non-real branches.\end{itemize}

\begin{theorem}\label{t3}
An $\AG$-diagram of an arbitrary real singularity determines its weak real topological type.
\end{theorem}

{\bf Proof.} Balke and Kaenders \cite{BK} proved that the $\AG$-diagram determines the topological type of
a totally real singularity, and we closely follow the lines of their proof referring for details to \cite[Section 2]{BK}
and presenting necessary modifications for the general case.

First, we remark that the partition requirement (see Section \ref{sec1}) was not, in fact, used in \cite{BK}. In particular, it is not needed
in the construction of the Coxeter-Dynkin diagram from the given divide as presented in \cite{GZ}.

\smallskip

{\bf(1)}
The main step in the proof of \cite[Theorem 2.5 and Corollary 2.6]{BK} is to show that
an $\AG$-diagram of a totally real singularity determines the branch structure of the divide, pairwise intersection numbers of the branches, and
an $\AG$-diagram of each branch. Their argument literally applies in the general case.
We notice in addition that one can easily distinguish between $\AG$-diagrams
of non-closed and closed branches of the divide, i.e., between an $\AG$-diagram of
a real branch of $(C,z)$ and an $\AG$-diagram of a pair of complex conjugate branches. Namely, in the former case,
the $\AG$-diagram contains either a univalent $\bullet$-vertex, or a bivalent $\bullet$-vertex joined with a $\oplus$-vertex and
$\ominus$-vertex, while in the latter case, the $\AG$-diagram has no such $\bullet$-vertices.

\begin{figure}
\setlength{\unitlength}{0.8cm}
\begin{picture}(14.5,7)(-1,-0.7)
\epsfxsize 135mm \epsfbox{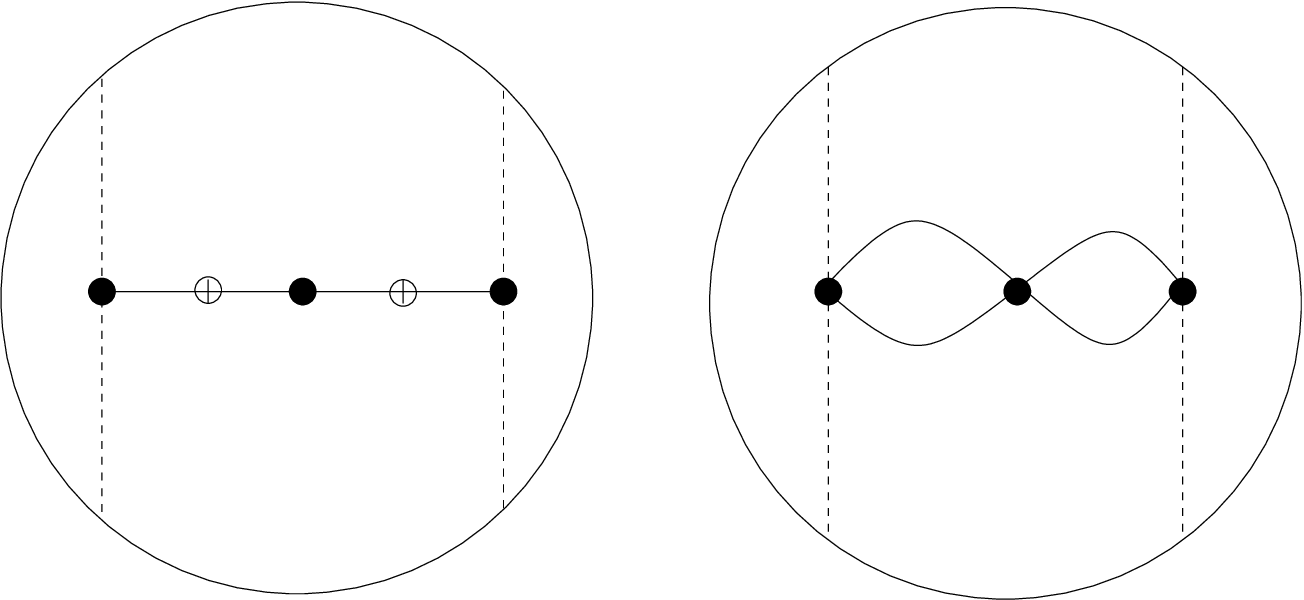}
\put(-13.2,-0.7){(a)}\put(-4,-0.7){(b)}\put(-7.2,3.8){$A$}\put(-0.8,3.8){$B$}
\put(-16.5,3.8){$A$}\put(-10,3.8){$B$}
\end{picture}
\caption{Chains of an $\AG$-diagram and of a divide} \label{fig2}
\end{figure}
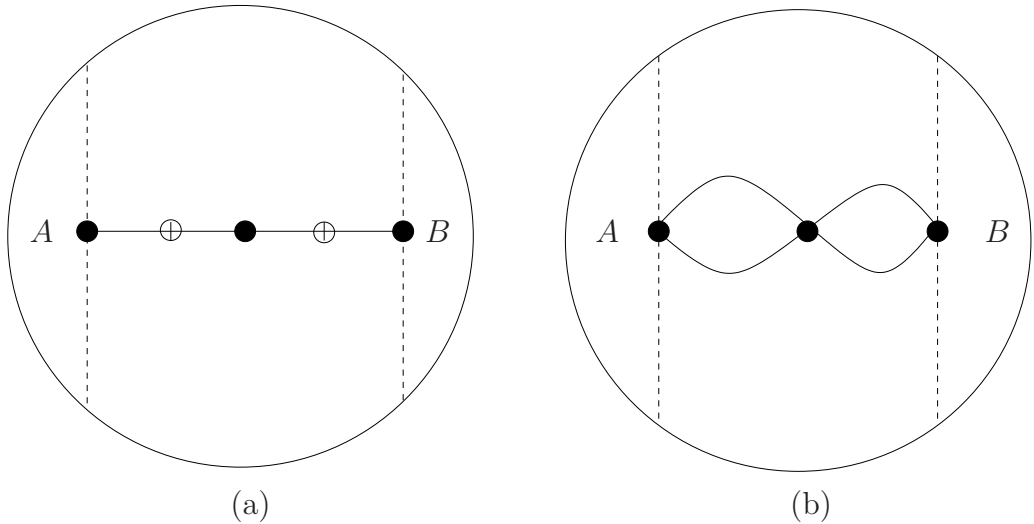

We only comment on the persistence of the cyclic order of real branches of the singularity (aka, non-closed branches of the divide).
An embedding of the $\AG$-diagram into $\R B_{C,z}$ defines the divide up to isotopy (see \cite[Page 46]{BK}).
The ambiguity in the construction of an embedding is related to the existence of the so-called
{\bf chains} in the $\AG$-diagram, i.e., connected subgraphs consisting of bivalent or univalent
$\bullet$-vertices and bivalent $\oplus$-vertices (or
bivalent $\ominus$-vertices) joined by arcs as shown in Figure \ref{fig2}(a) (cf. \cite[Figure 6]{BK}).
Figure \ref{fig2}(b) shows the corresponding fragment of the divide (cf. \cite[Figure 7]{BK}).
By \cite[Lemma 2.8]{BK}, the given $\AG$-diagram can be transformed by inserting new chains and extending the existing ones in a controlled way into a {\bf chain separating} $\AG$-diagram, whose maximal (with respect to inclusion) chains have pairwise distinct lengths, and no new chain can be added.

Each chain of a divide shares the boundary with two
non-inner components of the complement to the divide, and the disc $\R B_{C,z}$ can be cut into three parts as shown in
Figure \ref{fig2}(b) by dashed lines (cf. \cite[Figure 7]{BK}), and similarly one can cut
$\R B_{C,z}$ with respect to the embedded chain of the $\AG$-diagram, Figure \ref{fig2}(a). Then a given embedding of a chain separating $\AG$-diagram can be changed in part $A$ or in part $B$
by a reflection with respect to the axis of the chain
(and so for any other maximal chain).
Note that the branches of the divide, which are disjoint from the chain of the divide,
must all lie
either in part $A$, or in part $B$, since any two of them must intersect each other. In the presence of such branches, located, say, in part $A$,
and under the assumption that the chain is formed by two branches of the divide, all possible self-intersections of the
latter branches must lie in part $A$ too due to
Lemma \ref{l3}(i) applied to the divide with one of these two branches removed. All these observations yield that the cyclic order
of non-closed branches of the divide is preserved under the changes of the embedding of the chain separating $\AG$-diagram described above. Finally, we note that the same cycling order of the divide is induced by the corresponding embedding of the original $\AG$-diagram.

\smallskip

{\bf(2)} The topological type a real branch of the given singularity can be recovered from its $\AG$-diagram, see
\cite[Theorem 1.9]{BK}. In a similar way, we show that an $\AG$-diagram of a closed branch of the divide determines the topological type
of a real singularity formed by a pair of complex conjugate branches. Namely, an $\AG$-diagram defines the monodromy operator of such a singularity, see \cite{AC2} and
\cite[Page 39]{GZ1}, and hence its characteristic polynomial, which is the reduced Alexander polynomial of the link of the singularity
\cite[\S8]{M} (see also \cite[Theorem 3.3]{SW}). On the other hand, by \cite[Proposition 3.2]{CDGZ},
the reduced Alexander polynomial of a singularity formed by two topologically equivalent branches determines
the topological type of the branches and their intersection multiplicity.

\smallskip

{\bf(3)} To complete the recovery of the topological type of the given singularity
$(C,z)$, we have to
find pairwise intersection multiplicities of the branches of $(C,z)$. By \cite[Lemma 2.2]{BK}, the intersection number of two non-closed branches of the divide equals the intersection multiplicity of the corresponding real branches of $(C,z)$. Similarly, the intersection number of a non-closed and a closed branches of the divide equals twice the intersection multiplicity of the corresponding real branch of $(C,z)$ with each of the two complex conjugate branches of $(C,z)$ corresponding to the closed branch of the divide.
At last, consider the intersection of two closed branches of the divide and suppose without loss of generality that these are the only branches of the divide. From Lemma
\ref{l1203} we know the topological type and the intersection multiplicity of complex conjugate branches of $(C,z)$ associated with each of the branches of the divide. We claim that this information together with the intersection number of the branches of the divide determines the pairwise intersection multiplicities of all four branches of $(C,z)$. Indeed, this can easily be proved by induction on the number of real infinitely near points in the resolution tree of $(C,z)$.

\smallskip

{\small
\section*{Acknowledgements} The authors were supported by the grant no. 1174-197.6/2011
from the German-Israeli Foundations, by the grant no.
176/15 from the Israeli Science Foundation, and by a grant from the Hermann-Minkowski-Minerva
Center for Geometry at Tel Aviv University. The authors are indebted to Sergey Fomin
for prompting this research and for many stimulating discussions and important suggestions.
Special thanks are due to Ilya Tyomkin, who pointed out a gap in the first version of the paper.
We are also very thankful to the unknown referee for valuable remarks and corrections.}

\section*{Appendix}

We present here some material on Alexander polynomials of singularity links.

\medskip
{\bf (1)}
Given a singularity $(C,z)$ with $r$ local branches, its link $L=C\cap\partial B_{C,z}$ consists of $r$
algebraic knots in $\partial B_{C,z}\simeq S^3$ and it has a topological invariant
$\Delta^2_L(t_1,...,t_r)\in\Z[t_1,...,t_r]$ called the {\bf Alexander polynomial of the link} (see
\cite{SW} for precise definitions and detailed treatment). According to \cite[\S8, page 95]{M} (see also
\cite[Theorem 3.3]{SW}), the {\bf reduced Alexander polynomial}
$\Delta_L(t):=(t-1)\Delta^2_L(t,...,t)$ is the characteristic polynomial of the monodromy of $(C,z)$.

All the roots of the monic polynomial $\Delta_L(t)$ are known to lie on the unit circle, and hence, by Kronecker's theorem
\cite{Kro} (see also \cite[Section 34, Lemma (a)]{Hec}) are roots of unity. It follows that $\Delta_L(t)$ admits a (unique) factorization into the product of cyclotomic polynomials $\Phi_d(t)\in\Z[t]$\ \footnote{The roots of $\Phi_d$ are simple and they are all primitive roots of unity of degree $d$. The cyclotomic polynomials are distinct for distinct $d\ge1$ and are irreducible in $\Q[t]$.}.

Denote by $CT\subset\Q(t)^*$ the multiplicative subgroup generated by cyclotomic polynomials. For a function $f\in CT$, set $\MCTI(f)=d$ and $\MCTE(f)=\eps$, where $d$ is the maximal with the property that $\Phi_d(t)$ enters the product defining $f$, and $\eps$ is the exponent of $\Phi_d(t)$ in this product. Now, we construct a sequence of functions $f_i\in CT$, $i\ge1$, and integers as follows: Set $f_1(t)=\Delta_L(t)$, and for any $k\ge1$, inductively define
\begin{equation}d_k=\MCTI(f_k),\quad \eps_k=\MCTE(f_k),\quad f_{k+1}(t)=
f_k(t)(t^{d_k}-1)^{-\eps_k}\ .\label{e1703}\end{equation} Since the positive part of the sequence $\{d_k\}_{k\ge1}$ is strictly decreasing, there is the minimal $m(f)\ge1$ such that $f_n\equiv1$, $n\ge m(f)$.

\medskip

{\bf(2)} As we noticed in the end of Introduction, any A'Campo-Gusein-Zade diagram associated with a totally real singularity uniquely determines the (complex)
topological type of the singular point. The core of this phenomenon is the following statement due to W. Burau
\cite{Bu}, which we supply here with a proof.

\begin{lemma}\label{l1107}
The reduced Alexander polynomial of a unibranch singularity determines
the topological type of the singularity.
\end{lemma}

{\bf Proof.}
A real unibranch singularity is topologically equivalent to one given by the Puiseux series
$$y=x^{\frac{m_1}{n_1}}+...+x^{\frac{m_i}{n_1...n_i}}+
...+x^{\frac{m_s}{n_1...n_s}}\ ,$$
where the parameters are positive integers satisfying
$$s\ge1,\quad \gcd(m_j,n_j)=1\ \text{for all}\ j=1,...,s,\quad n:=n_1....n_s=\mt(C,z)\ ,
$$ \begin{equation}n_1\ge1,\ n_j>1\ \text{for all}\ j=2,...,s,\quad\text{and}\quad 1\le\frac{m_1}{n_1}<...<\frac{m_s}{n_1...n_s}\ .\label{e1107}
\end{equation}
The link $L:=C\cap\partial B_{C,z}$ is an
algebraic knot in $\partial B_{C,z}\simeq S^3$ and its Alexander polynomial can be written as
(see \cite[Section 8.3]{BrK})
\begin{equation}\Delta_L(t)=\frac{t-1}{t^n-1}\cdot
\prod_{j=1}^s\frac{t^{w_jb_{j,s}}-1}{t^{w_jb_{j+1,s}}-1}\ ,\label{e1407}\end{equation} where
\begin{equation}w_1=m_1,\quad w_j=m_j-m_{j-1}n_j+w_{j-1}n_jn_{j-1},\quad b_{j,s}=n_j...n_s,\quad 2\le j\le s\ .\label{ne1207b}\end{equation} If the singularity is a smooth curve germ, then
$\Delta_L(t)\equiv1$. Otherwise, $\deg\Delta_L(t)>0$ and we have
\begin{equation}1<n_1<m_1\ .\label{ne1207a}\end{equation} Observe that the inequalities in (\ref{e1107}) and (\ref{ne1207a}) imply
\begin{equation}n<w_1b_{2,s}\quad \text{and}\quad w_jb_{j+1,s}<w_jb_{j,s}<w_{j+1}b_{j+2,s}\ \text{for all}\ j\ge1,\label{e1207}\end{equation}
and then, in the associated with $\Delta_L(t)$ sequence $\{(f_k,d_k,\eps_k)\}_{k\ge1}$, we have
$$d_{2j-1}=w_{s+1-j}b_{s+1-j,s},\quad d_{2j}=w_{s+1-j}b_{s+2-j,s}, \quad j=1,...,s\ ,$$
$$d_{2s+1}=n,\quad d_{2s+2}=1\ .$$ It follows that $m(\Delta_L)=2s+3$, which yields the value of $s$. Furthermore, from the values of
$d_1,...,d_{2s}$ and formulas (\ref{ne1207b}), one can easily recover the parameters
$n_1,...,n_s$ and $m_1,...,m_s$, and thereby the topological type of the given singularity.
\proofend

\smallskip

In part (3) of the proof of Theorem \ref{t3} we used the following statement (which is, in fact, a particular case of \cite[Proposition 3.2]{CDGZ}). For the reader's convenience, we provide here a proof based on a simple direct computation.

\begin{lemma}\label{l1203} The reduced Alexander polynomial of a singularity formed by two topologically equivalent branches determines
the topological type of the branches and their intersection multiplicity.\end{lemma}

{\bf Proof.}
The singularity under consideration is topologically equivalent to a singularity $(C,z)$ with $z=0\in\C^2$ and two branches having the following Puiseux-type expansions:
$$y=x^{\frac{m_1}{n_1}}+...+x^{\frac{m_i}{n_1...n_i}}+\sqrt{-1}\left(x^{\frac{m_{i+1}}{n_1....n_{i+1}}}+
...+x^{\frac{m_s}{n_1...n_s}}\right)\ ,$$
$$y=x^{\frac{m_1}{n_1}}+...+x^{\frac{m_i}{n_1...n_i}}-\sqrt{-1}\left(x^{\frac{m_{i+1}}{n_1....n_{i+1}}}+
...+x^{\frac{m_s}{n_1...n_s}}\right)\ ,$$ where the parameters are positive integers satisfying
$$s\ge1,\quad 0\le i<s,\quad \gcd(m_j,n_j)=1\ \text{for all}\ j=1,...,s,\quad n:=n_1....n_s=\frac{1}{2}\mt(C,z)\ ,
$$ \begin{equation}n_1\ge1,\ n_j>1\ \text{for all}\ j=2,...,s,\ j\ne i+1,\quad\text{and}\ 1\le\frac{m_1}{n_1}<...<\frac{m_s}{n_1...n_s}\ .\label{e1203}
\end{equation}
Note that here all the pairs $(m_s,n_s)$, $s\ne i+1$, are characteristic Puiseux pairs; for $s=i+1$, the pair $(m_{i+1},n_{i+1})$
may be Puiseux pair as well and in this case, $n_{i+1}>1$, or may not and in this case $n_{i+1}=1$. This dichotomy reflects the position of the
last common infinitely near point of the two branches of the given singularity.

To recover the topological type of the branches and their intersection number, we need to know the
parameters
\begin{equation}s,\ i,\quad \text{and}\quad n_j,\ m_j,\ j=1,...,s\ .\label{e1603}\end{equation}

In our setting, the formula in \cite[Theorem 7.6]{SW} for the reduced Alexander polynomial of the link
$L=C\cap\partial B_{C,z}$ reads
\begin{equation}\Delta_L(t)=
\frac{t-1}{t^{2n}-1}\cdot\left(\prod_{j=1}^i\frac{t^{2w_jb_{j,s}}-1}{t^{2w_jb_{j+1,s}}-1}\right)
\cdot\frac{(t^{2w_{i+1}b_{i+1,s}}-1)^2}{t^{2w_{i+1}b_{i+2,s}}-1}
\cdot\left(\prod_{j=i+2}^s\frac{t^{n_je_j}-1}{t^{e_j}-1}\right)^2\ ,\label{e1803}\end{equation}
where $w_1,...,w_s$ are defined by (\ref{ne1207b}) and 
$$b_{j_1,j_2}=\prod_{j_1\le j\le j_2}n_j,\quad e_j=
w_{i+1}b_{i+1,s}b_{i+2,j-1}+w_jb_{j+1,s},\ i+2\le j\le s\ .$$

We can suppose that $(C,z)$ is not a node, since the node is easily recognized by the condition
$\deg\Delta_L(t)=\mu(C,z)=1$, and hence either $s>1$, or $m_1>1$. It follows from relations (\ref{e1203}) that
\begin{itemize}\item if $i\ge1$, then
\begin{equation}n<w_1b_{2,s}\quad \text{and}\quad w_jb_{j+1,s}<w_jb_{j,s}<w_{j+1}b_{j+2,s}\ \text{for all}\ 1\le j\le i,\label{e1303}\end{equation}
\item if $i+1<s$, then
\begin{equation}2w_{i+1}b_{i+1,s}<e_{i+2}\quad\text{and}\quad \begin{cases}e_j<n_je_j,\ &\text{for all}\ i+2\le j\le s,\\
n_je_j<e_{j+1}\ &\text{for all}\ i+2\le j<s,\end{cases}\label{e1403}\end{equation}
\item and
\begin{equation}w_{i+1}b_{i+2,s}\begin{cases}=w_{i+1}b_{i+1,s},\quad&\text{if}\ n_{i+1}=1,\\
<w_{i+1}b_{i+1,s},\quad&\text{if}\ n_{i+1}>1.\end{cases}\label{ne1507}\end{equation}
\end{itemize}
As in the proof of Lemma \ref{l1107}, we consider the auxiliary sequence $\{(f_k,d_k,\eps_k)\}_{k\ge1}$.
Observe, that, in the beginning, it has an even or odd number $l$ of even values of $\eps_k$ according as $n_{i+1}=1$, or $n_{i+1}>1$.
It follows that $s=[(m(\Delta_l)-1)/2]$ and $i+1=s-[l/2]$. Moreover, the sequence $d_k$, $k\ge1$, in (\ref{e1703}) provides values for
all the exponents of $t$ in the formula (\ref{e1803}). Considering this as a system of equations for $m_j,n_j$, $j=1,...,s$, we can easily resolve it and hence restore the topological type of the
given singularity.
\proofend

\medskip

{\bf(3)} The next question we address is: {\it How to recognize a polynomial which is the reduced Alexander polynomial of a real, irreducible over $\R$ singularity?}

\begin{lemma}\label{l1307}
(1) A monic polynomial $f(t)\in\Z[t]\cap CT$ of a positive degree is the Alexander polynomial of a unibranch singularity if and only if the associated sequence $\{(f_k,d_k,\eps_k)\}_{k\ge1}$ satisfies the following conditions:
\begin{itemize}\item $m(f)=2s+3$ with $s\ge1$;
\item $\eps_k=(-1)^{k-1}$ for all $k=1,...,2s$, $\eps_{2s+1}=-1$;
\item $d_{2i-1}/d_{2i}\in\Z$ and $\prod_{j=1}^id_{2j}/\prod_{j=1}^{i-1}d_{2j-1}\in\Z$
for all $i=1,...,s$, and $d_{2s+1}=\prod_{j=1}^s(d_{2j-1}/d_{2j})$.
\end{itemize}

(2) A monic polynomial $f(t)\in\Z[t]\cap CT$ of degree $>1$ is the reduced Alexander polynomial of a singularity consisting of two topologically equivalent local branches if and only if the associated sequence $\{(f_k,d_k,\eps_k)\}_{k=1,...,m}$ with $f_m\equiv1$ satisfies the following conditions:
\begin{itemize}\item either, $m(f)=2s+3$ for some $s\ge1$, and there exists $0\le i<s$ such that
$$\begin{cases}\eps_k=2(-1)^{k-1},\ &1\le k\le2s-2i-1,\\
\eps_k=(-1)^{k-1},\ &2s-2i\le k\le2s,\\\
\eps_{2s+1}=-1,\ &\eps_{2s+2}=1,\end{cases}$$
$$\frac{d_{2j-1}}{d_{2j}}\in\Z\quad \text{for}\ 1\le j\le s,\quad d_{2j}\cdot\prod_{k=1}^{j-1}\frac{d_{2k}}{d_{2k-1}}\in\Z\quad \text{for}\ 1\le j\le s-i-1,$$
$$\frac{1}{2}d_{2j}\cdot\prod_{k=1}^{j-1}\frac{d_{2k}}{d_{2k-1}}\in\Z\quad \text{for}\ s-i\le j\le s,\quad\text{and}\quad d_{2s+1}=2\prod_{j=1}^s\frac{d_{2j-1}}{d_{2j}}\ ;$$
\item or $m(f)=2s+2$ for some $s\ge1$, and there exists $0\le i<s$ such that
$$\begin{cases}\eps_k=2(-1)^{k-1},\ &1\le k\le2s-2i-2,\\
\eps_k=(-1)^k,\ &2s-2i\le k\le2s-1,\\\
\eps_{2s-2i-1}=1,\ &\eps_{2s}=-1,\quad\eps_{2s+1}=1,\end{cases}$$
$$\frac{d_{2j-1}}{d_{2j}}\in\Z\quad \text{for}\ 1\le 1\le s-i-1,\quad
\frac{d_{2j}}{d_{2j+1}}\in\Z\quad \text{for}\ s-i\le j\le s-1\ ,$$
$$\frac{1}{2}d_{2j+1}\cdot\prod_{k=1}^{s-i-1}\frac{d_{2k}}{d_{2k-1}}
\cdot\prod_{k=s-i}^{j-1}\frac{d_{2k+1}}{d_{2k}}\in\Z\quad\text{for}\ s-i-1\le j\le s-1\ ,$$
$$d_{2j}\cdot\prod_{k=1}^{j-1}\frac{d_{2k}}{d_{2k-1}}\in\Z\quad \text{for}\ 1\le j\le s-i-1,\quad\text{and}\quad d_{2s}=2\cdot\prod_{j=1}^{s-i-1}\frac{d_{2j-1}}{d_{2j}}\cdot
\prod_{j=s-i}^{s-1}\frac{d_{2j}}{d_{2j+1}}\ .$$
\end{itemize}
\end{lemma}

{\bf Proof.} {\bf(1)} Assuming that $f(t)$ is an Alexander polynomial, one immediately obtains from formula (\ref{e1407}) and relations (\ref{ne1207b})-(\ref{e1207}) all the conditions on $m(f)$, $\eps_k$, and $d_k$ mentioned in the hypotheses of part (1) of the lemma.

Assume now the hypotheses of part (1) of the lemma. Then
$$f(t)=\frac{t-1}{t^{d_{2s+1}}-1}\cdot\prod_{i=1}^s\frac{t^{d_{2i-1}}-1}{t^{d_{2i}}-1}\ .$$
Denote
$$n_i=\frac{d_{2s-2i+1}}{d_{2s-2i+2}},\quad w_i=\prod_{j=1}^{s-i+1}d_{2j}\left(\prod_{j=1}^{s-i}d_{2j-1}\right)^{-1},\quad i=1,...,s\ ,$$
which is equivalent to
$$d_{2j-1}=w_{s+1-j}\prod_{j=s+1-j}^sn_j,\quad d_{2j}=w_{s+1-j}\prod_{j=s+2-j}^s, \quad j=1,...,s\ .$$ In addition, we have $d_{2s+1}=n_1...n_s$. 
Now we recursively define
\begin{equation}m_i=m_{i-1}n_i+w_i-w_{i-1}n_{i-1}n_i,\quad i=1,...,s\ .\label{e1607}\end{equation} As required by (\ref{e1107}) and (\ref{ne1207a}), we have to show that
\begin{equation}n_1<m_1\quad\text{and}\quad m_i>m_{i-1}n_i,\ i=2,...,s\ .
\label{e1607a}\end{equation} Indeed, first,
$$d_{2s}>d_{2s+1}=\prod_{j=1}^s\frac{d_{2j-1}}{d_{2j}}\quad\Longrightarrow\quad n_1=\frac{d_{2s-1}}{d_{2s}}<m_1=w_1=\prod_{j=1}^sd_{2j}\left(\prod_{j=1}^{s-1}d_{2j-1}\right)^{-1}\ .$$ and, second,
$$m_i=m_{i-1}n_i+w_i-w_{i-1}n_{i-1}n_i$$
$$=m_{i-1}n_i+
\prod_{j=1}^{s-i+1}d_{2j}\left(\prod_{j=1}^{s-i}d_{2j-1}\right)^{-1}
-\prod_{j=1}^{s-i+2}d_{2j}\left(\prod_{j=1}^{s-i+1}d_{2j-1}\right)^{-1}
\cdot\frac{d_{3s-2i+3}}{d_{2s-2i+4}}\cdot\frac{d_{2s-2i+1}}{d_{2s-2i+2}}$$
$$=m_{i-1}n_i+(d_{2s-2i+2}-d_{3s-2i+3})\prod_{j=1}^{s-i}\frac{d_{2j}}{d_{2j-1}}>m_{i-1}n_i\ .$$

\medskip
{\bf(2)} Assuming that $f(t)$ is the reduced Alexander polynomial of a singularity consisting of two topologically equivalent local branches (and different from a node), from formula (\ref{e1803}) and relations (\ref{e1303})-(\ref{ne1507}), one can easily derive the conditions on $m(f)$, $\eps_k$, $d_k$ formulated in part (2) of the lemma.

In the opposite direction, given the conditions on $\eps_k$, $d_k$ corresponding to $m(f)=2s+3$, we obtain, first, that
$$f(t)=\frac{t-1}{t^{d_{2s+1}}-1}\cdot\left(\prod_{j=1}^i
\frac{t^{d_{2s-2j+1}}-1}{t^{d_{2s-2j+2}}-1}\right)
\cdot\frac{(t^{d_{2s-2i-1}}-1)^2}{t^{d_{2s-2i}}-1}
\cdot\left(\prod_{j=i+2}^s\frac{t^{d_{2s-2j+1}}-1}{t^{d_{2s-2j+2}}-1}\right)^2\ .
$$ Then we define the parameters
$$e_j=d_{2s-2j+2},\ i+2\le j\le s,\quad n_j=\frac{d_{2s-2j+1}}{d_{2s-2j+2}},\ 1\le j\le s,\quad b_{j_1,j_2}=\prod_{j_1\le k\le j_2}n_k\ ,$$
$$w_j=\frac{1}{2}d_{2(s-j+1)}\cdot(b_{j+1,s})^{-1},
\quad 1\le j\le i+1\ ,$$
$$w_j=d_{2(s-j+1)}\cdot(b_{j+1,s})^{-1}
-w_{i+1}b_{i+1,j}b_{i+2,j-1}\quad \text{for}\ i+2\le j\le s\ .$$ Since we are given $d_{2s+1}=2n_1...n_s$, the formula for $f(t)$ turns into the right-hand side of (\ref{e1803}).
So, it remains to consider the parameters $m_1,...,m_s$ given by (\ref{e1607}) and to prove inequalities (\ref{e1607a}). Up to $j\le i+1$, it follows by the argument used in the proof of part (1) of the lemma. If $j=i+2$, the required inequality amounts to
$$w_{i+1}n_{i+1}n_{i+2}<w_{i+2}$$
$$\Leftrightarrow\quad\frac{d_{2s-2i-1}}{2n_{i+1}...n_s}\cdot n_{i+1}n_{i+2}<
\frac{1}{n_{i+3}...n_s}\left(d_{2s-2i-2}-\frac{d_{2s-2i-1}}{2n_{i+1}...n_s}\cdot n_{i+1}...n_s\right)$$
$$\Leftrightarrow\quad d_{2s-2i-1}<d_{2s-2i-2}\ .$$ If $j>i+2$, the required inequality amounts to
$$w_{j-1}n_{j-1}n_j<w_j$$ which reads
$$\left(\frac{d_{2s-2j+4}}{n_j...n_s}-
w_{i+1}(n_{i+1}...n_{j-1})(n_{i+2}...n_{j-2})\right)n_{j-1}n_j$$ $$<
\frac{d_{2s-2j+2}}{n_{j+1}...n_s}-w_{i+1}(n_{i+1}...n_j)(n_{i+2}...n_{j-1})$$
$$\Leftrightarrow\quad d_{2s-2j+4}n_{j-1}=d_{2s-2j+3}<d_{2s-2j+2}\ .$$

Under the conditions on $\eps_k$, $d_k$ corresponding to $m(f)=2s+2$, we obtain that
$$f(t)=\frac{t-1}{t^{d_{2s}}-1}\cdot\left(\prod_{j=1}^i\frac{t^{d_{2s-2j}}-1}{t^{d_{2s-2j+1}}-1}\right)
\cdot(t^{d_{2s-2i-1}}-1)
\cdot\left(\prod_{j=i+2}^s\frac{t^{d_{2s-2j+1}}-1}{t^{d_{2s-2j+2}}-1}\right)^2\ .
$$ Then we define the parameters
$$n_j=\frac{d_{2s-2j}}{d_{2s-2j+1}},\ 1\le j\le i,
\quad n_{i+1}=1,\quad n_j=\frac{d_{2s-2j+1}}{d_{2s-2j+2}},\ i+2\le j\le s\ ,$$
$$b_{j_1,j_2}=\prod_{j_1\le k\le j_2}n_k,\quad e_j=d_{2s-2j+2},\ i+2\le j\le s\ ,$$
$$w_j=\frac{1}{2}d_{2(s-j+1)}\cdot(b_{j+1,s})^{-1},
\quad 1\le j\le i,\quad w_{i+1}=\frac{1}{2}d_{2s-2i-1}\cdot(b_{i+2,s})^{-1}\ ,$$
$$w_j=d_{2(s-j+1)}\cdot(b_{j+1,s})^{-1}
-w_{i+1}b_{i+2,j}b_{i+2,j-1},\quad i+2\le j\le s\ .$$
In the same manner as in the preceding paragraph, we get inequalities (\ref{e1107}) and (\ref{ne1207a}).
\proofend


\begin{thebibliography}{0}
\bibitem{AC} N. A'Campo.
Le groupe de monodromie des singularit\'es isol\'ees des courbes planes, I.
\textit{Math. Ann.} {\bf213} (1975), 1--32.
\bibitem{AC1} N. A'Campo. Real deformations and complex topology of plane curve singularities.
\textit{Annales de la Facult\'e de Science de Toulouse, 6-e Ser.,} {\bf 8} (1999), no. 1, 5--23.
\bibitem{AC4} N. A'Campo.
A combinatorial property of generic immersions of curves. \textit{Indag. Math. (N.S.)}
{\bf 11} (2000),  no. 3, 337--341.
\bibitem{AC2} N. A'Campo.
Quadratic vanishing cycles, reduction curves and reduction of the monodromy group of plane curve singularities.
\textit{Tohoku Math. J.} {\bf 53} (2001), 533--552.
\bibitem{AC3} N. A'Campo. Monodromy of real isolated singularities. \textit{Topology}
{\bf 42} (2003), no. 6, 1229--1240.
\bibitem{AGV} V. I. Arnold, A. N. Varchenko, and S. M. Gusein-Zade.
\textit{Singularities of Differentiable Maps, Vol. 2 (Monodromy and Asymptotics of Integrals)}.
Birkh\"auser, Boston, 1988.
\bibitem{BK} L. Balke and R. Kaenders. On a certain type of Coxeter-Dynkin diagrams of
plane curve singularities. \textit{Topology} {\bf 35} (1996), no. 1, 39--54.
\bibitem{Bra} K. Brauner. Zur Geometrie der Funktionen zweier Ver\"anderlicher.
\textit{Abh.
Math. Sem. Hamburg} {\bf 6} (1928), 8--54.
\bibitem{BrK} E. Brieskorn and H. Kn\"orrer. \textit{Plane algebraic curves}.
Birkh\"auser, Basel, 1986.
\bibitem{Bu} W. Burau. Kennzeichnung der Schlauchknoten.
\textit{Abh. Math. Sem. Hamburg} {\bf 9}, 125--133 (1933).
\bibitem{CDGZ} A. Campillo, F. Delgado, and S. M. Gusein-Zade.
On the topological type of a set of plane valuations with symmetries.
{\it Math. Nachr.} {\bf 290} (2017), no. 13, 1925--1938.
\bibitem{CP} O. Couture and B. Perron. Representative braids for links associated
to plane immersed curves. \textit{J. Knot Theory Ramifications} {\bf 9} (2000), no. 1, 1--30.
\bibitem{Gab} A. M. Gabrielov. Bifurcations, Dynkin diagrams, and modality of
isolated singularities. \textit{Func. Anal. Appl.} {\bf 8} (1974), no. 2, 94--98.
\bibitem{GLS} G.-M. Greuel, C. Lossen, and E. Shustin.
{\it Introduction to singularities and deformations}. Springer,
Berlin, 2007.
\bibitem{GZ2} S. M. Gusein-Zade. Dynkin diagrams for singularities of
functions of two variables. \textit{Func. Anal. Appl.} {\bf 8} (1974), no. 4, 295--300.
\bibitem{GZ} S. M. Gusein-Zade.
Intersection matrices for certain singularities of functions of two variables.
\textit{Func. Anal. Appl.} {\bf 8} (1974), no. 1, 10--13.
\bibitem{GZ1} S. M. Gusein-Zade. The monodromy groups of isolated singularities of hypersurfaces.
\textit{Russ Math. Surveys} {\bf 32} (1977), no. 2, 23--69.
\bibitem{Hec} E. Hecke. {\it Lectures on the Theory
of Algebraic Numbers}. Springer, 1981.
\bibitem{Kro} L. Kronecker. Zwei S\"atze \"uber Gleichungen mit ganzzahligen Coefficienten.
{\it J. reine und angew. Math.} {\bf53} (1857), 173--175.
\bibitem{M} J. Milnor. \textit{Singular points of complex hypersurfaces}.
Princeton Univ. Press, Princeton, 1968.
\bibitem{P} B. Perron. \textit{Preuve d'un th\'eor\`eme de N. A'Campo
sur les d\'eformations r\'eelles
des singularit\'es alg\'ebriques complexes planes}. Preprint,
Universit\'e de Bourgoigne, Dijon, 1998.
\bibitem{ShT} E. Shustin.   Gluing of singular and critical points.
\textit{Topology} {\bf 37} (1998), no. 1, 195--217.
\bibitem{Sh1} E. Shustin. A tropical approach to enumerative geometry.
\textit{St. Petersburg Math. J.} {\bf 17} (2006), 343--375.
\bibitem{ST} E. Shustin and I. Tyomkin. Patchworking singular algebraic curves, II.
\textit{Israel J. Math.} {\bf 151} (2006), 145--166.
\bibitem{SW} D. W. Sumners and J. M. Woods. The
Monodromy of Reducible Plane Curves. \textit{Invent. Math.} {\bf40} (1977), 107--141.
\bibitem{Z} O. Zariski. \textit{Algebraic Surfaces}, 2nd ed. Springer, Berlin etc., 1971.
\end{thebibliography}
\end{document}